\title{Contact Processes on Random Regular Graphs}
\author{
Steven Lalley and Wei Su\\
Department of Statistics\\
The University of Chicago\\
}
\date{\today}

\documentclass[10pt]{article}
\usepackage{mathrsfs}
\usepackage[letterpaper]{geometry}
\usepackage{amssymb, mathdots}
\usepackage{amsmath, turnstile}
\usepackage{amsfonts}
\usepackage{amsthm}
\usepackage{array}
\usepackage{gensymb}
\usepackage{multirow}
\usepackage{lscape}
\usepackage{color}
\usepackage{graphicx}
\usepackage{fullpage}
\usepackage{rotating}
\usepackage{listings}
\usepackage{float}
\usepackage{paralist}
\usepackage{caption}
\usepackage[stretch=10]{microtype}
\usepackage{tikz}
\usetikzlibrary[patterns]
\usetikzlibrary{arrows,shapes,trees, positioning} 
\usetikzlibrary{decorations.markings}
\tikzstyle{ann} = [fill=white,font=\footnotesize,inner sep=1pt]

\newtheorem{Theorem}{Theorem}[section]

\newtheorem{Lemma}[Theorem]{Lemma}

\newtheorem{Proposition}[Theorem]{Proposition}

\newcommand{\zz}[1]{\mathbb #1}
\newcommand{\goodevent}{F}

\begin{document}
\maketitle

\begin{abstract}
We show that the contact process on a random $d$-regular graph
initiated by a single infected vertex obeys the ``cutoff phenomenon''
in its supercritical phase. In particular, we prove that when the
infection rate is larger than the critical value of the contact
process on the infinite $d$-regular tree there are positive constants
$C, p$ depending on the infection rate such that for any $\varepsilon
>0$, when the number $n$ of vertices is large then (a) at times $t<
(C-\varepsilon)\log n$ the fraction of infected vertices is
vanishingly small, but (b) at time $(C+\varepsilon)\log n$ the
fraction of infected vertices is within $\varepsilon$ of $p$, with
probability $p$.

\end{abstract}

\section{Introduction}

The contact process with infection rate $\lambda>0$ on a connected,
locally finite graph $G=(\mathcal{V}_G,\mathcal{E}_G)$ is a continuous-time
Markov chain $(\xi_t)_{t\geq 0}$ with state space $\{\text{\rm subsets of }\mathcal{V}_G\}$
that evolves as follows:
(1) each infected site (that is, each vertex in $\xi_t$) recovers at
rate 1, and upon recovery is removed from $\xi_t$; and (2) each
healthy site (that is, each vertex not in $ \xi_t$) becomes infected
at rate $\lambda$ times the number of currently infected neighbors,
and upon infection is added to $\xi_t$. See \cite{liggett:book1} for a
formal construction, or alternatively \cite{harris} for the standard
\emph{graphical representation} using independent Poisson processes.

The behavior of the contact process on the infinite $d$-regular tree
$G=\zz{T}_{d}$ is reasonably well understood.  When $d=2$ (where
$\zz{T}_{2}=\zz{Z}$) there is a single survival phase
\cite{liggett:book1}, When $d\geq 3$, there are two survival phases:
in particular, there are critical values
$0<\lambda_1(\zz{T}_{d})<\lambda_2(\zz{T}_{d})<\infty$ such that
if $\lambda\leq \lambda_1$, then the contact process dies out almost
surely; if $\lambda_1<\lambda\leq\lambda_2$, then the contact process
survives globally with positive probability but dies out locally
almost surely; and if $\lambda>\lambda_2$ then the contact process
survives locally with positive probability. (See \cite{Pemantle} for
the case $d\geq 4$ and \cite{liggett2,stacey} for $d=3$.) The
parameter range $\lambda \in (\lambda_{1},\lambda_{2})$ is called the
\emph{weak survival phase}, and $\lambda>\lambda_{2}$ is the
\emph{strong survival phase}.

When $G$ is finite there is no survival phase, since the absorbing
state $\emptyset$ is accessible from every state $\xi \in
\{\text{subsets of } \mathcal{V}_G\}$.
Nevertheless, when the graph is large it will contain long
linear chains, and so if the infection rate is above the critical
value $\lambda_{1} (\zz{Z})$ the contact process will, with
non-negligible probability, survive for a long time in a
quasi-stationary state before ultimately dying out. This suggests
several problems of natural interest:

\medskip
\begin{compactenum}
\item [(A)] How does the survival time scale with the size of the graph?
\item [(B)] What is the nature of the quasi-stationary state?
\item [(C)] How does the process behave in its initial ``explosive'' stage?
\end{compactenum}

\medskip These questions have been studied for several important
families of graphs, notably the finite tori \cite{mountford}, finite
regular trees of large depth \cite{stacey2}, and versions of the
``small worlds'' networks of Watts and Strogatz
\cite{Durrett-Jung}. Stacey has shown that when $G_L$ is a finite
$d$-homogeneous rooted tree of depth $L$, the extinction time of a
contact process started from full occupancy in $G_L$ grows linearly in
$L$ when $\lambda<\lambda_2(\zz{T}_{d})$; but when
$\lambda>\lambda_2(\zz{T}_{d})$ it grows doubly exponentially in $L$,
and almost exponentially in the number of vertices. In a more recent
paper \cite{CMMV}, it has been proved that the extinction time grows
exponentially in the number of vertices.

In this paper we consider a different class of graphs, the
\emph{random $d$-regular graphs.} These are of interest for a number
of reasons: they are \emph{expanders}, they are locally tree-like, and
they are (unlike the finite trees) \emph{statistically homogeneous} in
a certain sense.  See \cite{wormald} for a survey. The behavior of
several common stochastic models on random $d-$regular graphs has been
studied. Lubetsky and Sly \cite{Lubetzky-Sly} have shown that the
simple random walk on a large random $d$-regular graph undergoes
\emph{cutoff}, that is, the transition to stationarity occurs in a
narrow time window.  Bhamidi, Hofstad, and Hooghiemstra
\cite{Bhamidi-Hofstad-Hooghiemstra}
have shown that distance between two randomly chosen vertices in first
passage percolation on a random $d-$regular graph is concentrated
around $C\log n$. Chatterjee and Durrett \cite{Chatterjee-Durrett}
have shown that the \emph{threshold} contact process on a random
$d-$regular graph exhibits a phase transition in the infectivity
parameter. More recently, Ding, Sly, and Sun \cite{Ding-Sly-Sun} have shown that
the independence number of a random d-regular graph is sharply
concentrated about its median.

The main result of this paper is that the contact process on a random
regular graph exhibits a cutoff phenomenon analogous to that for the
simple random walk.  We shall assume throughout that $nd$ is even and
$d\geq 3$.  Let $G\sim \mathcal{G} (n,d)$ be a random graph uniformly
distributed over the set of all $d$-regular graphs on the vertex set
$V=[n]$.  Given $G$, for any subset $A\subset [n]$, let $\xi^{A}_{t}$
be a contact process run on $G$ with initial state $\xi^{A}_{0}=A$
(when $A=\{u\}$ is a singleton, we will write $\xi^{u}_{t}$ for the
contact process with initial state $\xi^{u}_{0}=\{u \}$).  We shall be
primarily interested in the ``meta-stable'' phase, where the infection
rate $\lambda$ exceeds $\lambda_1(\zz{T}_{d})$, and our main focus
will be the following question: for a typical pair of vertices, what
is the time needed for a contact process started from one vertex to
infect the other? Since the diameter of a typical random regular graph
is on the order of $\log n$, we expect the infection time to be on the
same order. The main result of this paper, Theorem \ref{THM1}, implies
that this conjecture is true.

We say that a property holds \emph{asymptotically almost surely} if the
set of graphs in $\mathcal {G}(n,d)$ satisfying this property has
probability approaching 1 as $n$ goes to infinity. To denote
conditional probabilities and expectations given the graph $G$,
we will use a subscript $\zz{P}_{G}$ and $\mathbb{E}_G$.
Our main results are the following two theorems.

\begin{Theorem}\label{THM1}
Fix arbitrary $u\neq v \in [n]$, and let
$G\sim\mathcal {G}(n,d)$ be a random $d$-regular graph on the vertex set
$V=[n]$.  For each $\lambda >\lambda_{1} (\zz{T}_{d})$ let
$p_{\lambda}>0$ and $c_{\lambda}>0$  be the survival probability and
growth rate (cf. Proposition \ref{CP.CPExpectationgrowthrate}) for the
contact process started from the root on the infinite tree
$\zz{T}_{d}$ with infection rate
$\lambda$. Then for any  $0<\varepsilon<1/8$ there exist constants $g_{n}
(\varepsilon)\rightarrow 0$ as $n \rightarrow \infty$ such that for
asymptotically almost every $G$,
\begin{equation}\label{eq:mainResultLower}
\zz{P}_G\{ v\in \xi^u_s \quad \text{for some} \;\; s\leq
(1-\varepsilon)c_{\lambda}^{-1}\log n \} \leq g_n (\varepsilon)
\end{equation}
and
\begin{equation}\label{eq:mainResultUpper}
	\zz{P}_G\{v\in \xi^u_{(1+\varepsilon)c_{\lambda}^{-1}\log
	n}\}\geq(1-g_{n} (\varepsilon))p^{2}_{\lambda }.
\end{equation}
\end{Theorem}

\begin{Theorem}\label{THM2}
Let $\xi^G_{t}$ be a contact process with initial configuration
$[n]$. Fix $\varepsilon >0$; then for each $\delta>0$ there exist
constants $f_n(\delta)\rightarrow 0$ as $n\rightarrow \infty$ such
that for asymptotically almost every $G$,
\[\zz{P}_G\{ (1-\delta)np_{\lambda}\leq |\xi^G_{(1+\varepsilon)c_\lambda^{-1}\log n}|\leq (1+\delta)np_{\lambda} \}\geq 1-f_n(\delta),\]
\end{Theorem}

Assertion \eqref{eq:mainResultLower} will be proved in section 3, and
assertion \eqref{eq:mainResultUpper} in section 4. Theorem \ref{THM2}
will be proved in section 5. Throughout this paper we fix $\lambda>\lambda_1(\zz{T}_d)$.

\section{Preliminaries: Contact Process on the Infinite Regular Tree}
In this section, $\xi_t$ will denote the contact process started from
a single vertex $O$ (the \emph{root}) on the $d$-regular tree
$\zz{T}_{d}$.  The $d-$regular tree is a non-amenable graph, in the
sense that its Cheeger constant is positive. This can be expressed as
follows.  For a finite subset $S$ of vertices of $\zz{T}_{d}$, call
$v\in S$ a \emph{border point} if among the $d$ connected components
obtained by removing $v$ from $\zz{T}_{d}$, at least one of them
contains no other vertices in $S$. Let $B(S)$ be the set of border
points in $S$; then
\[|B(S)|\geq \left(1-\frac{1}{d-1}\right)|S|.\]
See, for instance,  Lemma 6.2 of \cite{Pemantle} for a proof.
We will denote by $h(\zz{T}_{d})$  the constant $1-1/(d-1)$.

The nonamenability of $\zz{T}_{d}$ implies that the supercritical
contact process on $\zz{T}_{d}$ grows exponentially. Here is a precise
formulation, proved in \cite{Madras-Schinazi} and
\cite{Morrow-Schinazi-Zhang}.

\begin{Proposition}\label{CP.CPExpectationgrowthrate}
There exist constants $c_{\lambda}>0$ and $C_d>0$ such that
\[\exp{(c_{\lambda}t)}\leq \zz{E}|\xi_t|\leq C_d\exp{(c_{\lambda}t)}.\]
\end{Proposition}

We will make frequent use of an auxiliary process, the \emph{severed}
contact process. We follow the terminology and notation of
\cite{Morrow-Schinazi-Zhang} and \cite{Pemantle}.  Define a branch
$\mathscr{B}$ to be the connected component of the root in the
subgraph obtained by removing a distinguished subset of $d-1$ edges,
each having an endpoint at the root $O$. The severed contact process
is the contact process restricted to $\mathscr{B}$, i.e., infection is
not allowed to travel across any of the $d-1$ removed edges. We will
use the letter $\eta$ to denote the severed contact process; in
particular, $\eta^S_t$ is the severed contact process with initial
configuration $S\subset \mathscr{B}$, and
$\eta_t=\eta^{O}_t=\eta^{\{O\}}_t$ the severed contact process
started with $O$ infected at time 0. Clearly, the severed contact
process $\eta_t$ is stochastically dominated by the contact process
$\xi_t$. In the standard graphical representation \cite{harris}
contact process $\xi_t$ and the severed contact process are
naturally coupled in such a way that $\eta^{S}_t\subset \xi^{S}_t$
for any initial configuration $S$ and all $t\geq 0$. Hence the
expected cardinality of infected sites in the severed contact process
is no larger than that of the original one. However, the following
proposition states that the severed process has comparable cardinality
in expectation.

\begin{Proposition}\label{CP.severedgrowthrate}
There exists a constant $A_1=A_1(\lambda,d)>0$  such that
\[\zz{E}|\eta_t|\geq A_1\exp{(c_\lambda t)}, \forall t\geq0.\]
\end{Proposition}

\begin{proof}
Inequality (5) of \cite{Morrow-Schinazi-Zhang} states, in our notation, that
\[\int_0^t \zz{E}|\eta_s|ds\geq \frac{1}{\lambda} \left( \frac{1}{d} \zz{E}|\xi_t|-1\right).\]
On the other hand, it is easy to see that
\[\zz{E}|\eta_s| \leq \zz{E}|\xi_s|\leq \frac{C_d}{d}\exp{(c_{\lambda}s)}.\]
Now fix a constant $T>0$ (to be determined later); we have
\begin{equation*}
\begin{aligned}
\frac{1}{\lambda}\left( \frac{1}{d} \zz{E}|\xi_{nT}|-1\right)&\leq \int_0^{nT} \zz{E}|\eta_s|ds\\
&\leq\int_0^{(n-1)T}\frac{C_d}{d}\exp{(c_{\lambda}s)}ds+\int_{(n-1)T}^{nT}\zz{E}|\eta_s|ds\\
&=\frac{C_d}{dc_{\lambda}}\left(\exp{(c_\lambda(n-1)T)}-1\right)+\int_{(n-1)T}^{nT}\zz{E}|\eta_s|ds,\\
\end{aligned}
\end{equation*}
and therefore
\begin{equation*}
\begin{aligned}
\int_{(n-1)T}^{nT}\zz{E}|\eta_s|ds&\geq \frac{1}{d\lambda} \zz{E}|\xi_{nT}|-\frac{1}{\lambda}-\frac{C_d}{dc_{\lambda}}\left(\exp{(c_\lambda(n-1)T)}-1\right)\\
&= \exp{(c_\lambda nT)}\left(\frac{1}{d\lambda}-\frac{C_d}{dc_{\lambda}}\frac{1}{\exp{(c_\lambda T)}}\right)+\left(\frac{C_d}{dc_{\lambda}}-\frac{1}{\lambda}\right).
\end{aligned}
\end{equation*}
This inequality holds for any $T>0$, so we can choose $T$ large enough
that
$\frac{1}{d\lambda}-\frac{C_d}{dc_{\lambda}}\frac{1}{\exp{(c_\lambda
T)}}>0$.  Fix such  $T$; then there exists $t_n^{*}\in
[(n-1)T,nT]$ such that
\[\zz{E}|\eta_{t_n^*}|\geq \frac{1}{T} \left(\exp{(c_\lambda
nT)}\left(\frac{1}{d\lambda}-\frac{C_d}{dc_{\lambda}}\frac{1}{\exp{(c_\lambda
T)}}\right)+\left(\frac{C_d}{d}-\frac{1}{\lambda}\right)\right).\]

Now we use the fact
\[\zz{E}|\eta_{t+s}|\geq (h(\zz{T}_{d})\zz{E}|\eta_s|-1)\zz{E}|\eta_t|.\]
This holds because if we run the severed contact process up to time
$s$, and then keep only those infected vertices that are border points
of $\eta_s$, and run severed contact processes from each of these
border points for another $t$ inside the unoccupied branch at time
$s$, the resulting infection set is dominated by the original severed
contact process at time $t+s$. (The extra $-1$ is because the origin
$O$ has some edges removed and we might not be able to run a severed
contact process from it) Therefore, for any time $t\in [nT,(n+1)T]$,
we have
\[\zz{E}|\eta_t|\geq (h(\zz{T}_{d}) \zz{E}|\eta_{t_n^*}|-1)
\zz{E}|\eta_{t-t_n^*}|\geq (h(\zz{T}_{d}) \zz{E}|\eta_{t_n^*}|-1)
\times\inf_{s\in [\,0,2T]}\zz{E}|\eta_s|,\]
and so we conclude that
\[\inf_{t\geq 0} \frac{\zz{E}|\eta_t|}{\exp{(c_{\lambda}t)}} >0.\]
\end{proof}

Next we will show that, conditioned on survival up to a large time
$t$, the cardinality of the contact process is concentrated
around its expectation, at least in exponential rate.

\begin{Proposition}\label{CP.CPgrowthrate}
For any fixed $\delta>0$, we have
\begin{equation}\label{eq:aboveAndBelow}
\zz{P}\{\exp{(c_\lambda(1-\delta)t)}\leq |\xi_t|\leq
\exp{(c_\lambda(1+\delta)t)} \ \vert  |\xi_t|>0\}\rightarrow
1,\text { \rm  as } t\rightarrow \infty.
\end{equation}
\end{Proposition}

The proof relies on another useful estimate, Theorem 4 of Athreya
\cite{Athreya}. Athreya's result is proved for the case where $p_0=0$ (where $(p_i)_{i\in
\mathbb{N}}$ is the offspring distribution of the branching process),
but it can be  generalized painlessly to the case $p_0>0$. We record
this extension of Athreya's theorem as  Lemma \ref{CP.Athreya}

\begin{Lemma}\label{CP.Athreya}
Suppose $(Z_n)_{n\in \mathbb{N}}$ is a branching process with mean
offspring number $\mu>1$, and such that  $\zz{E}\left(\exp(\theta_0
Z_1)\ \vert Z_0=1\right)<\infty$ for some $\theta_0>0$. Then there
exists  $\theta_1>0$ such that
\[\sup_{n\geq 1} \zz{E} \exp\left(\theta_1 \frac{Z_n}{\mu^n}\right)<\infty.\]
\end{Lemma}

\begin{proof}[Proof of Proposition \ref{CP.CPgrowthrate}]
We will deal with the two inequalities in the event in \eqref{eq:aboveAndBelow} separately.

(A) First, we will show that
\begin{equation}\label{eq:aboveAndBelowUnconditional}
\zz{P}\{|\xi_t|\leq \exp{(c_\lambda(1+\delta)t)}
\}\rightarrow 1, \text{ \rm as } t\rightarrow \infty.
\end{equation}
Since the event $\{|\xi_t|>0\}$ has probability bounded away from $0$,
it will follow that the conditional probability given $\{|\xi_t|>0\}$
also converges to $1$.

To prove \eqref{eq:aboveAndBelowUnconditional}, we will build a
discrete-time branching random walk on the tree which stochastically
dominates the contact process.  Specifically, we define
$(BRW_{nT})_{n\in\mathbb{N}}$ as follows:

(1) $BRW_0=\xi_0=\{O\}$ and  $BRW_T=\xi_T$ (the state of the  contact process at time $T$).

(2) Given $BRW_{nT}$, for each particle in $BRW_{nT}$, run independent
contact processes for time $T$ starting at the locations of these
particles. Notice
that we allow multiple particles to occupy the same vertex. The set of
all particles (and their locations) is defined to be $BRW_{(n+1)T}$.

The process $BRW_{nT}$ is a discrete-time branching random walk, and
so its cardinality $(|BRW_{nT}|)_n$ is a Galton-Watson branching
process. The mean offspring number is $\zz{E}|\xi_T|\leq
C_d\exp(c_\lambda T)$. Moreover, the distribution of $|BRW_T|$ is
dominated by a geometric distribution, because the contact process is
dominated by a Yule process, which has a geometric distribution at any
specific time, so the finite moment generating function assumption in
Lemma \ref{CP.Athreya} holds here. Consequently, Lemma
\ref{CP.Athreya} implies that if $T$ is sufficiently large, then
\begin{equation*}
\zz{P}\{ |\xi_{nT}|\geq
\exp(c_\lambda(1+\delta)nT)\}\rightarrow 0, \text{ \rm as } n
\rightarrow \infty.
\end{equation*}
Since the branching random walk stochastically dominates the contact
process, the result \eqref{eq:aboveAndBelowUnconditional} follows, at
least for $t$ in the arithmetic progression $\{nT \}_{n\geq 0}$. To
extend \eqref{eq:aboveAndBelowUnconditional} to all $t$, we use a
simple ``filling'' argument, as follows.

Suppose that there exists a sequence of time points $t_n\rightarrow
\infty$ such that $\zz{P}\left( |\xi_{t_n}|\geq
\exp(c_\lambda(1+\delta)t_n)\right)\geq \varepsilon>0$ for some
$\varepsilon$, and without loss of generality $t_n\in
[k_nT,(k_n+1)T)$. Assume that $T$ is large enough so that not
only \[\zz{P}\{ |\xi_{nT}|\geq
\exp(c_\lambda(1+\delta)nT)\}\rightarrow 0, \text{ \rm as } n
\rightarrow 0,\] but also
\[\zz{P}\{ |\xi_{nT}|\geq \exp(c_\lambda(1+\delta/2)nT)\}\rightarrow 0,\text{ \rm as } n \rightarrow 0.\]
Given that $|\xi_{t_n}|\geq \exp(c_\lambda(1+\delta)t_n)$, for each
infection at $t_n$, there is  a positive (and fixed) probability
$p_T$ that it remains alive for at least time $T$, and therefore there
is a constant order lower bound for the probability that at
time $(k_n+1)T$, there are at least
$0.99p_T\exp(c_\lambda(1+\delta)k_nT )$ infections. This would
contradict the fact that
\[\zz{P}\{ |\xi_{nT}|\geq \exp(c_\lambda(1+\delta/2)nT)\}\rightarrow 0, \text{\rm  as } n \rightarrow 0,\]
because when $n$ is large enough, $0.99p_T\exp(c_\lambda(1+\delta)nT
)>\exp(c_\lambda(1+\delta/2)(n+1)T)$.

(B) Next, we show that
\begin{equation}\label{CP.1}
\zz{P}\{|\xi_t|\geq \exp{(c_\lambda(1-\delta)t)} \ \vert |\xi_t|>0\}\rightarrow 1,\text{\rm as }t\rightarrow \infty
\end{equation}
by showing that for each $\delta >0$ there exists $T>0$ such that
\begin{equation}\label{CP.2}
\zz{P}\{\liminf_{n\rightarrow \infty}
\log(|\xi_{nT}|)/(nT) \geq C_\lambda(1-\delta) \ \vert |\xi_t|>0,
\forall t>0 \}= 1.
\end{equation}
Together with another ``filling'' argument, this will prove that
(\ref{CP.2}) holds along the entire time axis, and (\ref{CP.1}) will
follow easily.  We will prove (\ref{CP.2}) by establishing the
following two assertions:

\begin{compactenum}
\item [\textsc{Assertion} 1:] On the event that a severed contact process survives,
the exponential growth rate is as desired.
\item [\textsc{Assertion} 2:] On the event that a contact process survives, we can
find a surviving severed contact process embedded in it.
\end{compactenum}

For Assertion 1, recall that for a severed contact process $\eta_t$, we
have  $\zz{E}|\eta_t|\geq A_1 \exp(c_\lambda t)$, by Proposition
\ref{CP.severedgrowthrate}. Now  construct a branching process as
follows:

\begin{compactenum}
\item [(1)] At time 0 there is  1 particle at the root.
\item [(2)] Run the severed contact process up to time $T$, then keep only
the infections on the ``border", i.e., infections which  have at least
one uninfected branch (at time $T$) of the tree connected to it.
\item [(3)] For all remaining infections, specify one uninfected branch for
each infection, and run a severed contact process within this branch
with the infection serving as the new ``root". Repeat (2) and (3).
\end{compactenum}
Denote the cardinality of infected vertices at time $nT$ as
$\{X_n\}$. The above branching process has mean offspring number

\[\zz{E}X_1\geq h(\zz{T}_{d}) \zz{E}|\eta_T|-1\geq h(\zz{T}_{d})A_1
\exp(c_\lambda T)-1\geq \exp(c_\lambda(1-\delta)T),\]
provided $T$ is large enough. It is clear that $X_n\leq |\eta_{nT}|\leq |\xi_{nT}|$.
By the Kesten-Stigum theorem, \cite{K-S},  if a branching
process $(X_n)_{n\in\mathbb{N}}$ with mean offspring number $\mu$
satisfies $\zz{E} X_1\log_{+}X_1<\infty$, then $\lim
X_n/\mu^n>0$ almost surely on the event of survival. This proves Assertion 1.

Assertion 2 is even easier. From Assertion 1, we know that the chance that
the severed contact process survives is $p>0$, so on the event
$|\xi_t|\geq M$, there will be at least $\lfloor
h(\zz{T}_{d})M\rfloor$ border points at time $t$, and the chance there
is at least one surviving severed contact process is at least
$1-(1-p)^{\lfloor h(\zz{T}_{d})M\rfloor}$. On the event that the
contact process survives, we know $|\xi_t|\rightarrow \infty$ as
$t\rightarrow \infty$, so the chance that we can find a surviving
severed contact process is 1.
\end{proof}

As we will see in later sections, we will explore the random regular
graph with the growth of the contact process, and it is important to
know up to time $t$ how big is the explored set compared to the
cardinality of the infected set at time $t$. Therefore we would like
to investigate the growth rate of the quantity $|\cup_{s\leq t}
\xi_s|$. The following proposition states that up to a constant
factor, it is comparable to $|\xi_t|$.

\begin{Proposition}\label{CP.Unionexpectationgrowthrate}
There exists $B_1=B_1(\lambda,d)<\infty$ such that
\[\exp(c_\lambda t)\leq \zz{E}|\cup_{s\leq t}\xi_s|\leq B_1 \exp(c_\lambda t), \forall t\geq 0.\]
\end{Proposition}

\begin{proof}
The first inequality $\exp(c_\lambda t)\leq \zz{E}|\cup_{s\leq
t}\xi_s|$ follows directly from Proposition
\ref{CP.CPExpectationgrowthrate}. For the second inequality,
we claim that
\begin{equation}\label{CP.5}
\zz{E}|\cup_{u\leq t+s}\xi_u|\leq \zz{E}|\cup_{u\leq
t}\xi_u|+\zz{E}|\xi_t|\times \zz{E}|\cup_{u\leq s}\xi_u|.
\end{equation}
This holds because if a site has been infected by time $t+s$, there
are two possibilities: either it was infected by time $t$, or it is
infected during the time interval $[t,t+s]$ by an infection alive at
time $t$. These account for the two terms on the right hand side of
(\ref{CP.5}).

Thus, for any fixed $T>0$, we have
\[
\zz{E}|\cup_{u\leq (n+1)T}\xi_u|\leq \zz{E}|\cup_{u\leq nT}\xi_u|+\zz{E}|\xi_{nT}|\times \zz{E}|\cup_{u\leq T}\xi_u|,
\]
and
\begin{equation*}
\begin{aligned}
&\frac{\zz{E}|\cup_{u\leq (n+1)T}\xi_u|}{\zz{E}|\xi_{(n+1)T}|}\\
\leq &\frac{\zz{E}|\cup_{u\leq nT}\xi_u|+\zz{E}|\xi_{nT}|\times \zz{E}|\cup_{u\leq T}\xi_u|}{\zz{E}|\xi_{(n+1)T}|}\\
\leq &\frac{\zz{E}|\cup_{u\leq nT}\xi_u|+\zz{E}|\xi_{nT}|\times \zz{E}|\cup_{u\leq T}\xi_u|}{1/C_d^2\zz{E}|\xi_{nT}|\zz{E}|\xi_{T}|} \quad \text{from Proposition \ref{CP.CPExpectationgrowthrate}},\\
=&\frac{C_d^2}{\zz{E}|\xi_T|}\frac{\zz{E}|\cup_{u\leq nT}\xi_u|}{\zz{E}|\xi_{nT}|}+\frac{C_d^2\zz{E}|\cup_{u\leq T}\xi_u|}{\zz{E}|\xi_T|},
\end{aligned}
\end{equation*}
so if we let $r_n=\zz{E}|\cup_{u\leq nT}\xi_u|/\zz{E}|\xi_{nT}|$, then
\[r_{n+1}\leq \frac{C_d^2}{\zz{E}|\xi_T|}r_n+\text{constant}.\]
As long as $C_d^2/\zz{E}|\xi_T|<1$, which holds if we take $T$ large
enough, we conclude that $\sup_n r_n<\infty$, that is,
\[\sup_n \frac{\zz{E}|\cup_{u\leq nT}\xi_u|}{\zz{E}|\xi_{nT}|}<\infty.\]
Using an argument similar to that  in the proof of Proposition
\ref{CP.severedgrowthrate}, we can extend this to
\[\sup_{t\geq 0}\frac{\zz{E}|\cup_{s\leq t}\xi_s|}{\zz{E}|\xi_t|}<\infty,\]
which translates to the desired inequality, by Proposition
\ref{CP.CPExpectationgrowthrate}.
\end{proof}

\begin{Proposition}\label{CP.Uniongrowthrate}
For any fixed $\delta>0$, we have
\[\zz{P}\{\exp{(c_\lambda(1-\delta)t)}\leq |\cup_{s\leq t}\xi_s|\leq \exp{(c_\lambda(1+\delta)t)} \ \vert  |\xi_t|>0\}\rightarrow 1, \text{ \rm as } t\rightarrow \infty.\]
\end{Proposition}

\begin{proof}
This follows by the same argument as in
proof of part (i) of Proposition \ref{CP.CPgrowthrate}.
\end{proof}

\begin{Proposition}\label{CP.SeveredUniongrowthrate}
There exists $B_1=B_1(\lambda,d)<\infty$ such that
\[A_1\exp(c_\lambda t)\leq \zz{E}|\cup_{s\leq t}\eta_s|\leq B_1 \exp(c_\lambda t), \forall t\geq 0,\]
where $A_1$ is from Proposition \ref{CP.severedgrowthrate}.
\end{Proposition}

\begin{proof}
This is an immediate corollary of Proposition
\ref{CP.Unionexpectationgrowthrate}.
\end{proof}

Next we will introduce an important concept that will figure
prominently in the arguments of sections 4 and 5. For any vertex $x\in
\xi_t$, say that $x$ is a \emph{pioneer point} if $x\in B(\cup_{s\leq
t}\xi_s)$, in other words, there exists a branch of the tree connected
to $x$ which has been completely uninfected up to time $t$. We call
such a branch a \emph{free branch}.  Notice that a pioneer point is
automatically a border point of $\xi_t$. We will use $\zeta_t$ to
denote the collection of pioneer points at time $t$. The next
proposition describes the approximate size of $\zeta_t$ conditional on
survival.

\begin{Proposition}\label{CP.pioneerpoint}
For any $\delta>0$,
\[\mathbb{P}\{| \zeta_t|\geq \exp(c_\lambda (1-\delta)t)   \, \vert \,|\xi_t|>0\}
\rightarrow 1, \text{\rm as }t\rightarrow \infty.\]
\end{Proposition}
\begin{proof}
From Proposition \ref{CP.CPgrowthrate} and  \ref{CP.Uniongrowthrate},
we may assume that conditional on the event $\{ |\xi_t|>0 \}$,
both of the following events occur:
\begin{compactenum}
\item [(1)] $|\xi_t|\geq  \exp(c_\lambda (1-\delta) t)$;
\item [(2)] $|\cup_{s\leq t} \xi_s| \leq \exp(c_\lambda (1+\delta)t)$.
\end{compactenum}
Assuming these, it is easy to deduce that there are at
least $(1-o(1))\exp(c_\lambda(1-\delta)t)$ vertices in $\xi_t$ such
that for each such vertex there is a branch connected to it which contains
no more than $\exp(3c_\lambda \delta t)$ vertices in $\cup_{s\leq t}\xi_s$.
Fix such a vertex $x\in\xi_t$ and such a branch that has no more than
$\exp(3c_\lambda \delta t)$ vertices in $\cup_{s\leq t}\xi_s$. Notice
that the vertices that belong to $\cup_{s\leq t}\xi_s$ in this branch are
connected. Therefore, there exists a path $y_0 y_1\dots y_L$ in this branch,
where $L\leq \log_{d-1}(\exp(3c_\lambda \delta t)) \leq 6\delta c_\lambda t$,
such that
\begin{compactenum}
\item [(1)] $y_0=x$;
\item [(2)] $y_{i}$ is connected to $y_{i+1}$ in this branch, for
$0\leq i\leq L-1$;
\item [(3)] $y_L$ is connected to a branch that has no vertices in $\cup_{s\leq t}\xi_s$.
\end{compactenum}
Now run the contact process for another time $6\delta c_\lambda t$,
and hope that at the end there is some chance of creating a pioneer
point at $y_L$ by infections along the path $y_0y_1\dots y_L$. We will
say that we have a \emph{successful infection event} if all of the
following events happen:
\begin{compactenum}
\item [(1)] between time $[i,i+1]$, the infection at $y_i$ infects $y_{i+1}$
before itself dies, and after getting infected, the vertex $y_{i+1}$ neither
dies nor infects other vertices in the interval $[i,i+1]$
\item [(2)] after $y_L$ first becomes gets infected between time $[L-1,L]$, this infection
neither dies, nor infects other vertices before time $6\delta c_\lambda t$.
\end{compactenum}
The (conditional) probabilities of (1) and (2) are bounded away from $0$.
Given that a  successful infection event occurs, the vertex  $y_L$
 becomes a pioneer point in $\zeta_{(1+6\delta c_\lambda )t}$.
On the other hand, a successful infection event will happen with chance at least
$q^{6\delta c_\lambda t+1}$, where $q=\min(q_1,q_2)>0$.

For different $y_0$'s, their corresponding successful infection events
are mutually independent, and so the number of successful infection events
stochastically dominates  the  Binomial distribution
\[W\sim \text{Binomial} (\, (1-o(1))\exp(c_\lambda(1-\delta)t), q^{6\delta c_\lambda t+1}), \]
and as long as $\delta$ is sufficiently small, with probability approaching 1,
\[W> (1-o(1)) \exp(c_\lambda(1-\delta)t) q^{6\delta c_\lambda t+1}>\exp(c_\lambda (1-D\delta )t),\]
for some $D>0$.
Therefore,  conditional on the event $|\xi_t|>0$,
with probability approaching 1,
\[|\zeta_{(1+6\delta c_\lambda)t}|>\exp(c_\lambda(1-D\delta )t). \]
This is essentially the desired conclusion if  $\delta$ is
small enough.
\end{proof}

Similarly, for the severed contact process $\eta_t$, we can define
pioneer points to be those vertices in $\eta_t$ that are also border
points of $\cup_{s\leq t}\eta_s$.  Denote the
set of such pioneer points by $\psi_t$.
By the same proof as for Proposition
\ref{CP.pioneerpoint}, we obtain the following proposition.

\begin{Proposition} \label{CP.severedpioneerpoint}
For any $\delta>0$,
\[\mathbb{P}\{| \psi_t|\geq \exp(c_\lambda (1-\delta)t)  \,\ \vert |\eta_t|>0\}
\rightarrow 1, \text{\rm as }t\rightarrow \infty.\]
\end{Proposition}

Finally, we show that the event that the severed contact process grows
exponentially faster than it is supposed to is exponentially
unlikely.

\begin{Proposition}\label{CP.exponentialbound}
For all $\delta>0$, there exists a constant $K>0$ and $\gamma>1$ such that
\[\zz{P}\{ |\cup_{s\leq t}\eta_s|\geq \exp((1+\delta)c_\lambda
t)\}\leq \exp(-K \gamma^t), \text{ \rm as } t\rightarrow
\infty.\]
\end{Proposition}
\begin{proof}
This is an immediate result from Proposition \ref{CP.Athreya} and
\ref{CP.SeveredUniongrowthrate}, once we observe that $(|\cup_{s\leq
nT}\eta_s|)_n$ is dominated by a branching process with mean offspring
number $\zz{E}|\cup_{s\leq T}\eta_s|$.
\end{proof}

\section{Contact Process on a Random Regular Graph}

\subsection{The configuration model for random regular graphs} A
\emph{random $d$-regular graph} is a graph chosen uniformly from the
collection $\mathcal{G} (n,d)$ of all $d$-regular graphs on the vertex
set $[n]$. We assume that $dn$ is even.  A useful way to construct a
random $d$-regular graph is the \emph{configuration model} introduced
by Bollob\'as \cite{Bollobas} (also see \cite{Bollobas:book} and
\cite{wormald}).  This works as follows. To each of the $n$ vertices
$u$, associate $d$ distinct half-edges $(u,i)$, and perform a uniform
perfect matching on these $dn$ half-edges. Using this matching,
construct a (multi-)graph by placing an edge between vertices $u$ and
$v$ for every pair of half-edges $(u,i)$ and $(v,j)$ that are matched.
The resulting graph need not be connected, and it might have multiple
edges and self-loops; however, the probability that the configuration
model produces a simple, connected graph is bounded away from 0 as $n
\rightarrow \infty$ (cf. \cite{wormald}). Moreover, given that the
resulting graph is simple (that is, has no self-loops or multiple
edges), it is uniformly distributed over $\mathcal {G}(n,d)$. Thus,
whenever an event holds w.h.p. for the (multi-)graph obtained from the
configuration model, it also holds w.h.p. under the uniform
distribution on $\mathcal {G}(n,d)$.

An important feature of the configuration model is that, at any stage,
the first half-edge in the next random pair can be selected using
any rule, as long as the second half-edge is chosen uniformly at
random from the remaining half-edges (see \cite{wormald}).

\subsection{Growth and exploration constructions}

\subsubsection{Vanilla version}

There are two layers of randomness in our model: first, the graph $G$
is chosen from the uniform distribution on the set $\mathcal{G} (n,d)$
of $d$-regular graphs with vertex set $[n]$, and then a contact
process is run on $G$. We would like to construct a probability space
where we use the configuration model to grow the contact
process and the random graph $G$ in tandem, building edges of $G$ only
at those times $t$ when the contact process attempts a new infection
from a vertex whose neighborhood structure is not yet completely
determined.

Suppose we want to run a contact process $(\xi^u_t)_t$ with initial configuration $\{u\}$.
Let $U_t$ be the set of unmatched half-edges
up to time $t$. At time $0$, only vertex $u$ is infected, and no edges are
yet determined, so $U_0$ is full, that is, it contains all $nd$ half-edges.
The recovery and infection
times of the contact process are determined by a system of independent
Poisson processes attached to the vertices of the graph, two to each
vertex (one for recoveries, the other for outgoing infections).

At any time $t$ when an infected vertex $v$ attempts an infection, one of the
$d$ half-edges incident to $v$ is selected uniformly at random. If this half-edge
is already matched to a half-edge $(w,j)$ then vertex $w$ is infected, if it
is currently healthy, or left infected if already infected. If, on the
other hand, one of the unmatched half-edges $(v,i)$ incident to $v$ is
selected then one of the other remaining unmatched half-edges $(w,j)$ is
chosen at random from $U_t\backslash \{(v,i)\}$ and matched with $(v,i)$,
and vertex $w$ is infected.  After this we remove $(v,i)$ and $(w,j)$ from $U_t$.

We will refer to this construction as the \emph{vanilla version} of the
 \emph{grow and
explore} process, and denote by $\zz{P}$ and $\zz{E}$ the probability
and expectation operator for the probability space on which the
underlying Poisson processes and other random variables used in the
matchings and infection attempts are defined. If we run the grow
and explore process up to time $T$, it is possible that
$U_T$ will not be not empty; in this case we match the remaining half-edges in $U_T$
to complete the graph  $G$. We condition on the event that the obtained graph is simple.

\begin{Proposition}\label{growandexplore.equivalence}
If the grow and explore process is run up to time $T$, then,
conditional on the event that the resulting graph $G$ is simple, the
pair ($G, (\xi_t)_{0\leq t \leq T})$ will have the same joint
distribution as for the contact process on a random regular graph.
\end{Proposition}

\begin{proof}
First of all, $G$ is uniform over $\mathcal{G}(n,d)$, because
whenever we pair two unmatched half-edges, the second half-edge is always
chosen uniformly at random from the unmatched pool.

Secondly, the interoccurence times between infection attempts are
i.i.d. exponentials, and in each attempt the active infection chooses
one of its neighbors at random and independent of everything else.
Meanwhile each vertex allows at most one infection at a time. Therefore
$(\xi_t)_{0\leq t\leq T}$ is a version of the contact process on $G$.
\end{proof}

The above construction of the grow and explore process assumes a
singleton initial configuration and that the graph is initially
completely unexplored. It is clear that the construction can be
trivially modified so as to work with an arbitrary initial
configuration and with part of the graph initially explored.

\subsubsection{Cover tree version with singleton initial configuration}

Next we describe a variation of the grow and explore process, in which
we grow a contact process $\tilde{\xi}_{t}$ on the infinite cover tree
and in tandem assign labels $v\in [n]$ to the vertices of
$\zz{T}=\zz{T}_{d}$ in such a way that $\tilde{\xi}$ partially
projects, via the labeling $\phi $ of vertices, to a contact process $\xi_{t}$
on $G$. The assignment of labels to vertices of $\zz{T}_{d}$ will result in a
(random) \emph{labeling function}
\begin{align*}
	\phi :\mathcal{V}_{\zz{T}_{d}} \longrightarrow [n] 
\end{align*}
that will determine the covering map from $\zz{T}_{d}$ to $G$ \emph{and}
the edge structure of $G$, as well as the projection mapping from the
vertex set of $\zz{T}_{d}$ to that of $G$. The construction will require
that some -- but not all -- of the vertices in $\tilde{\xi }_{t}$ be
selected for projection to $\xi_{t}$; thus, at any time $t$ the set
$\tilde{\xi}_{t}$ will be partitioned as $\tilde{\xi}_t
=\tilde{\xi}_{t,\text{\textrm{BLUE}}}\cup \tilde{\xi}_{t,\text{RED}}$, and
\[
	\xi_{t}=\phi (\tilde{\xi}_{t,\text{\textrm{BLUE}}}).
\]
Where appropriate, we will denote vertices of $\zz{T}_{d}$ via a tilde,
e.g., $\tilde{x}$, and use $x$ to denote the corresponding vertex
$x\in [n]$, so that $\phi (\tilde{x})=x$.

Fix a vertex $u\in [n]$; the singleton $\{u\}$ will be the initial
configuration of the (projected) contact process on $G$. Denote the
root vertex of the infinite tree $\zz{T}_{d}$ by $\tilde{u}$, and declare
$\phi (\tilde{u})=u$. Let $\tilde{\xi}_{t}$ be a contact process on
$\zz{T}_{d}$ with initial configuration $\tilde{\xi}_{0}=\tilde{u}$;
assume that this is constructed in the usual way, using independent
Poisson processes attached to the vertices $\tilde{v}$ of $\zz{T}_{d}$ to
determine the times at which recoveries and attempted infections
occur, and independent uniform random variables to determine which
neighbor of a vertex $\tilde{v}$ will be selected when vertex
$\tilde{v}$ attempts an infection. At time $t=0$ only the label $\phi
(\tilde{u})=u$ is determined; the function $\phi $ is augmented only
at those times when a \emph{blue} vertex of $\tilde{\xi}_{t}$ attempts
an infection.  The rules by which these augmentations occur are as
follows.

Suppose at time $t$, an infected vertex $\tilde{x}\in \tilde{\xi
}_{t,\textrm{BLUE}}$ attempts an infection. At this time $t$ some of the
neighbors of $\tilde{x}$ might have been labeled, and others might
not; denote by $\tilde{y}_1,\dots, \tilde{y}_\ell$ the neighbors that
have been labeled, with $\phi (\tilde{y}_{i})=y_{i}$, and by
$\tilde{z}_1,\tilde{z}_2,\dotsc ,\tilde{z}_{d-\ell}$ the neighbors that have
not yet been labeled. Moreover, at this time some of the neighbors of
$x$ in $G$ will have been determined, \emph{including}
$y_{1},y_{2},\dotsc ,y_{\ell}$, but possibly also some others, which
we denote by $y_{\ell +1},y_{\ell+2},\dotsc ,y_{\ell +k}$, where $\ell
+k \leq d$. Because the infection attempt entails choosing one of the
$d$ neighbors of $\tilde{x}$ at random to serve as the target of the
attempt, there are 3 possibilities:

\begin{compactenum}
\item [(1)] With probability $l/d$,  one of the vertices $\tilde{y}_1,\dots,\tilde{y}_\ell$ is
chosen. In this case $\phi$ is not augmented.
\item [(2)] With probability $k/d$, one of the vertices  $\tilde{z}_1,\dots,\tilde{z}_{d-\ell}$
is chosen randomly, say $\tilde{w}$, and one of the labels $y_{\ell+1},\dots,
y_{\ell+k}$ is chosen uniformly  at random to serve as the
label $\phi (\tilde{w})$ for the vertex $\tilde{w}$.
\item [(3)] With probability $1-l/d-k/d$, an unused half-edge
 $(x,i)$ incident to $x$ is chosen randomly, and another unused
 half-edge $(w,j)$ then chosen randomly from among all remaining
 unused half-edges  and matched with $(x,i)$.
 Then one of the vertices $\tilde{z}_1,\dots,\tilde{z}_{d-\ell}$'s is
 randomly selected and  labeled $w$. When this happens  we
add an edge connecting  $x$ and $w$ to $G$ and remove the two half-edges $(x,i)$ and $(w,j)$
from the set of unused half-edges.
\end{compactenum}

To complete the construction, we must specify how the vertices of
$\tilde{\xi }_{t}$ are to be colored (red or blue). This is done as
follows. First, a vertex $\tilde{v}$ is assigned a color only at a
time $t$ when it enters (or re-enters, if it was previously infected
but subsequently recovered) $\tilde{\xi }_{t}$. Second, if $\tilde{v}$
is infected by a \emph{red} vertex $\tilde{x}$ then it is colored
\emph{red}. Third, if $\tilde{v}$ is infected by a \emph{blue} vertex
$\tilde{x}$ then it is colored \emph{blue} unless the label $v$
assigned to $\tilde{v}$ has also been assigned to one of the other
vertices of $\tilde{\xi}_{t,\textrm{BLUE}}$; in this case $\tilde{v}$
is colored \emph{red}.

Notice that in this construction, when $\phi (\tilde{\xi}_{t,\text{\textrm{BLUE}}})$ attempts
an infection, any of its $d$ neighbors is equally likely to be the
target of the infection attempt, and when a new edge
is added to $G$  it follows the configuration model. Therefore, the
projection obeys the same rules as in the vanilla version of grow and
explore described above. This proves the following.

\begin{Proposition}\label{growandexplore.equivalence2}
The pair $(G, (\phi(\tilde{\xi}_{t,\text{\textrm{BLUE}}}) )_{0\leq t\leq T})$
obtained by running the cover tree version of grow and explore
has the same law as in the vanilla version of the grow and
explore process.
\end{Proposition}

The cover tree version of grow and explore has  two constituent processes:
the contact process $\tilde{\xi}_t$ on the cover tree and
the labeling process. We will call these the \emph{grow} process
and the \emph{explore} process, respectively.

To emphasize the initial configuration $\{u\}$, its corresponding
contact process on the cover tree is denoted as $\tilde{\xi}^u_t$.
Later we will run several contact processes on multiple vertices,
and adding the superscript will help us distinguish them.

The following figures illustrate concrete examples of how the construction works.

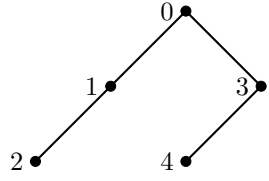
\begin{figure}[H]
\centering
\begin{tikzpicture}[scale=0.5]
\node[draw, circle, scale=0.4, fill=black] (u) at (0,0) {};
\node[draw, circle, scale=0.4, fill=black] (u) at (-2,-2) {};
\node[draw, circle, scale=0.4, fill=black] (u) at (-4,-4) {};
\node[draw, circle, scale=0.4, fill=black] (u) at (2,-2) {};
\node[draw, circle, scale=0.4, fill=black] (u) at (0,-4) {};
\draw[thick](0,0) -- (-4,-4)--cycle;
\draw[thick](0,0) -- (2,-2)--cycle;
\draw[thick](2,-2) -- (0,-4)--cycle;
\node[]  at (-0.5, 0) {$0$};
\node[]  at (-2.5, -2) {$1$};
\node[]  at (-4.5, -4) {$2$};
\node[]  at (1.5, -2) {$3$};
\node[]  at (-0.5, -4) {$4$};
\end{tikzpicture}
\caption {This graph records the order that vertices on a 3-regular cover tree are first infected. $0$ indicates the vertex is in the initial configuration.}\label{fig:ex1}
\end{figure}

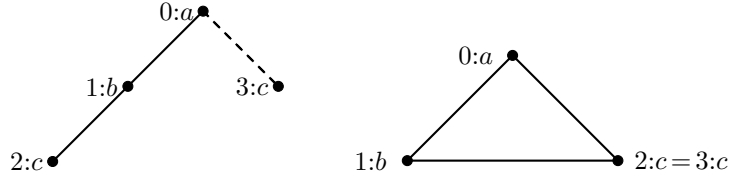
\begin{figure}[H]
\centering
\begin{tikzpicture}[scale=0.5]
\node[draw, circle, scale=0.4, fill=black] (u) at (0,0) {};
\node[draw, circle, scale=0.4, fill=black] (u) at (-2,-2) {};
\node[draw, circle, scale=0.4, fill=black] (u) at (-4,-4) {};
\node[draw, circle, scale=0.4, fill=black] (u) at (2,-2) {};
\draw[thick](0,0) -- (-4,-4)--cycle;
\draw[dashed, thick](0,0) -- (2,-2)--cycle;
\node[]  at (-0.7,0) {0:$a$};
\node[]  at (-2.7,-2) {1:$b$};
\node[]  at (-4.7,-4) {2:$c$};
\node[]  at (1.3,-2) {3:$c$};
\end{tikzpicture}
\qquad
\begin{tikzpicture}[scale=0.7]
\node[draw, circle, scale=0.4, fill=black] (u) at (0,0) {};
\node[draw, circle, scale=0.4, fill=black] (u) at (-2,-2) {};
\node[draw, circle, scale=0.4, fill=black] (u) at (2,-2) {};
\draw[thick](0,0) -- (-2,-2)--cycle;
\draw[thick](0,0) -- (2,-2)--cycle;
\draw[thick](-2,-2) -- (2,-2)--cycle;
\node[]  at (-0.7,0) {0:$a$};
\node[]  at (-2.7,-2) {1:$b$};
\node[]  at (3.2,-2) {2:$c$\,=\,3:$c$};
\end{tikzpicture}
\caption {Given Figure \ref{fig:ex1}, we try to project it onto the finite 3-regular graph. Suppose the top vertex has already been labelled as $a$, we sequentially label vertex 1,2,3,4 according to the law we describe before. When we label vertex 3 we accidentally use label $c$ again, so the finite graph immediately becomes the one on the right. Remember that we only allow one infection per vertex, so whenever we observe multiple infections at the same time on a vertex we will remove the whole infection trail except the chronically first one.}\label{fig:ex2}
\end{figure}

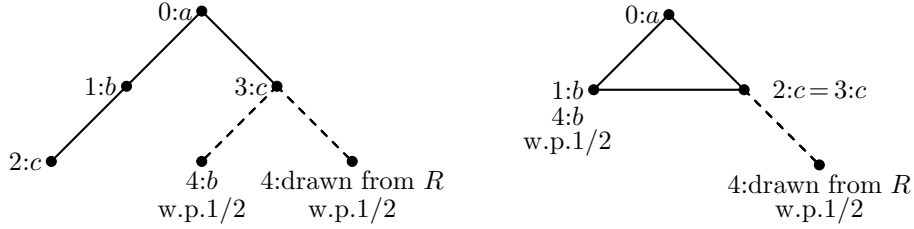
\begin{figure}[H]
\centering
\begin{tikzpicture}[scale=0.5]
\node[draw, circle, scale=0.4, fill=black] (u) at (0,0) {};
\node[draw, circle, scale=0.4, fill=black] (u) at (-2,-2) {};
\node[draw, circle, scale=0.4, fill=black] (u) at (-4,-4) {};
\node[draw, circle, scale=0.4, fill=black] (u) at (2,-2) {};
\node[draw, circle, scale=0.4, fill=black] (u) at (0,-4) {};
\node[draw, circle, scale=0.4, fill=black] (u) at (4,-4) {};
\draw[thick](0,0) -- (-4,-4)--cycle;
\draw[thick](0,0) -- (2,-2)--cycle;
\draw[dashed, thick](2,-2) -- (0,-4)--cycle;
\draw[dashed, thick](2,-2) -- (4,-4)--cycle;
\node[]  at (-0.7,0) {0:$a$};
\node[]  at (-2.7,-2) {1:$b$};
\node[]  at (-4.7,-4) {2:$c$};
\node[]  at (1.3,-2) {3:$c$};
\node[]  at (0,-4.5) {4:$b$ };
\node[]  at (0,-5.3) {w.p.1/2};
\node[]  at (4,-4.5) {4:drawn from $R$ };
\node[]  at (4,-5.3) {w.p.1/2};
\end{tikzpicture}
\qquad
\begin{tikzpicture}[scale=0.5]
\node[draw, circle, scale=0.4, fill=black] (u) at (0,0) {};
\node[draw, circle, scale=0.4, fill=black] (u) at (-2,-2) {};
\node[draw, circle, scale=0.4, fill=black] (u) at (2,-2) {};
\draw[thick](0,0) -- (-2,-2)--cycle;
\draw[thick](0,0) -- (2,-2)--cycle;
\draw[thick](-2,-2) -- (2,-2)--cycle;
\node[]  at (-0.7,0) {0:$a$};
\node[]  at (-2.7,-2) {1:$b$};
\node[]  at (4,-2) {2:$c$\,=\,3:$c$};
\node[draw, circle, scale=0.4, fill=black] (u) at (4,-4) {};
\draw[dashed,thick](2,-2) -- (4,-4)--cycle;
\node[]  at (-2.7,-2.7) {4:$b$};
\node[]  at (-2.7,-3.4) { w.p.1/2 };
\node[]  at (4,-4.5) {4:drawn from $R$  };
\node[]  at (4,-5.2) {w.p.1/2 };

\end{tikzpicture}
\caption {Given Figure \ref{fig:ex1} and \ref{fig:ex2}, suppose the infection trail coming from vertex 3 to vertex 4 is not removed, then when we label vertex 4, according to our law it has chance 1/2 to be $b$ and chance 1/2 to be drawn from $R$. The left graph is on the cover tree and the right one is on the finite graph.}\label{fig:ex3}
\end{figure}

\subsection{ Proof of Assertion (\ref{eq:mainResultLower}) of Theorem \ref{THM1}}

Let $\xi_t^u=\phi(\tilde{\xi}_{t,\textrm{BLUE}})$ be the contact process
constructed from the cover tree version of the grow and explore
process.  Let $t_1=(1-\varepsilon)\log n/c_\lambda$.  To prove
\eqref{eq:mainResultLower}, it suffices to show that
\begin{equation}\label{eq:proofLower}
	\mathbb{P} \{\exists s\leq t_1 \,\text{ \rm such that }\,
	v\in\xi^u_s \}\rightarrow 0, \text{\rm as } n\rightarrow
	\infty.,
\end{equation}
because if $\mathbb{P} \{\exists s\leq t_1 \,s.t.\,
\xi^u_s(v)=1 \}\leq l_n$ and $l_n\rightarrow 0$ as
$n\rightarrow\infty$, then by Markov's inequality,
\[
\mathbb {P}\{G: \zz{P}_G(\exists s\leq t_1 \,\text{\rm such that }\,
v\in \xi^u_s )\geq \sqrt{l_n}\}\leq \sqrt{l_n}.
\]
Now the event that there exists $s\leq t_{1}$ such that $v\in
\xi^{u}_{s}$ coincides with the event that at least one vertex in
$\cup_{s\leq t_1}\xi^u_{s}$ is labelled $v$. Consequently, to bound
the probability of this event,  it suffices to show that

\begin{compactenum}
\item [(1)]  as $n\rightarrow\infty$, $\zz{P}\{|\cup_{s\leq
t_1}\tilde{\xi}^u_{s}|\leq n^{1-\,\varepsilon/2}\}\rightarrow 1$;
and
\item [(2)] Given $|\cup_{s\leq t_1}\tilde{\xi}^u_{s}|\leq
n^{1-\,\varepsilon/2}$, the conditional probability that label $v$ is not
used in the labeling process before time $t_{1}$ approaches $1$.
\end{compactenum}

Assertion (1) follows directly from Proposition
\ref{CP.Uniongrowthrate}. To prove (2), observe that, because all
labels other than $u$ are equally likely to be used in the labeling
process, and because these probabilities add up to at most
$n^{1-\,\varepsilon/2}$, the chance that label $v$ appears in
$\cup_{s\leq t_1}\tilde{\xi}^u_{s}$ is at most
$n^{1-\,\varepsilon/2}/(n-1)\rightarrow 0$ as $n\rightarrow \infty$.
\qed

\subsection{Cover tree version with general initial configuration}

In the next section we will build the contact process with generic
initial configuration, therefore we need to generalize our grow and
explore process. Suppose the initial configuration of the contact process
is a set $\xi_0=\{u_1,u_2,\dots,u_k\}\subset [n]$, and some of the
half-edges are already paired up. We will construct a cover tree for each
of the $u_i$'s, and assign a labeling function at each of these cover trees,
i.e., we have
\[\phi_1\times \phi_2\times\dots \times \phi_k : \mathcal{V}_{\zz{T}_{d}}^k \longrightarrow
[n]^k,\]
where $\mathcal{V}_{\zz{T}_{d}}^k$ is the cartesian product of $k$ copies of
$\mathcal{V}_{\zz{T}_{d}}$ and $[n]^k$ is the cartesian product of $k$ copies of $[n]$.

At each $u_i$, we associate an independent contact process on $\zz{T}_{d}$.
We call it $(\tilde{\xi}_t^{u_i})_{t\geq 0}$, with $\tilde{\xi}_0^{u_i} =\tilde {u}_i$
and $\phi_i (\tilde{u}_i)=u_i$.
Similar to the definition in section 3.2, each $\tilde{\xi}_t^i$ is partitioned as
$\tilde{\xi}_{t,\text{BLUE}}^i \cup \tilde{\xi}_{t,\text{RED}}^i$, and
\[ \xi_t =\bigcup_{i=1}^k \phi_i (\tilde{\xi}_{t,\text{BLUE}}^i). \]

When there is an infection attempt in $\tilde{\xi}_{t,\text{BLUE}}^i$,
the law of updating the labeling function $\phi_i$ is exactly the same
as in section 3.2, such that whenever $\xi_t$ attempts an infection,
any of its $d$ neighbors is equally likely to be the target. Notice
now all $\phi_i$'s share the same unused pool of half-edges.
The rule for coloring  $\tilde{\xi}_t^i$ is essentially the same as in
section 3.2, i.e., descendants of red vertices are red;
descendants of blue vertices are blue unless there is another vertex
in $\tilde{\xi}_{t,\text{BLUE}}^j$ for some $j$ that has the same
label assigned to it, in which case the newborn infection is colored
red.  It is not hard to obtain the following analogy of Proposition
\ref{growandexplore.equivalence2}.

\begin{Proposition}\label{growandexplore.equivalence3}
The pair $(G, (\xi_t )_{0\leq t\leq T})$ obtained by
running the above version of grow and explore
has the same law as in the vanilla version of the grow and
explore process.
\end{Proposition}
In following sections we slightly abuse notations by ignoring
the index of $\phi_i$ and use $\phi$ as the only labeling function.

Similarly we can substitute the process on the tree with the severed
contact process and use the same labeling law. This will be useful in
the next section.

\subsection{Building independent contact processes}

Suppose now we would like to build several independent contact processes,
for example two independent contact processes with different initial configuration,
$(\xi^u_t)_{0\leq t \leq T_1}$ and $(\xi^v_t)_{0\leq t \leq T_2}$ on the same
random regular graph. To do so, we can first run the grow and explore
process to obtain $(\tilde{\xi}^u_t)_{0\leq t \leq T_1}$ and the labeling
function $\phi$. After that, we run an independent contact process on
the cover tree, $(\tilde{\xi}^v_t)_{0\leq t \leq T_2}$. In order to construct
the labeling function for $(\tilde{\xi}^v_t)_{0\leq t \leq T_2}$, we have to
be consistent with what has been explored on the finite graph so far.
That is, besides obeying the rule of constructing the labeling function in
previous sections, the pool of unused half-edges is whatever remained
after labeling $(\tilde{\xi}^u_t)_{0\leq t\leq T}$. In other words the processes
on the cover trees are independent, while the labeling process have to
be mutually consistent.

It is not hard to see that the processes obtained in this manner have
the desired distribution. Moreover, the order of labeling the processes
on the tree does not matter: we can instead label $(\tilde{\xi}^v_t)_{0\leq t\leq T_2}$
and then $(\tilde{\xi}^u_t)_{0\leq t\leq T_1}$ without changing the
distribution of the processes on the finite graph.

\section{A Second Moment Argument}

\subsection{Heuristics and Strategy}\label{ssec:strategy} In this section we shall
prove assertion~\eqref{eq:mainResultUpper} of Theorem~\ref{THM1}. This
states that for any two vertices $u,v\in [n]$ the conditional
probability, given the graph $G$, that $v\in \xi^{u}_{t_{+}}$
converges to $p_{\lambda}^{2}$ as $n \rightarrow \infty$, where
$t_{+}:= (1+\varepsilon)c_{\lambda}^{-1}\log n$. The rationale, in brief,
is as follows. For most vertices $u$ and most $d$-regular graphs $G$
the contact process $\xi^{u}_{t}$ on $G$ looks -- at least locally --
like a contact process on the infinite regular tree
$\zz{T_{d}}$ initiated by a single infected vertex $\tilde{u}$
at the root. The chance that such a contact process survives is
$p_{\lambda}$, so the chance that the contact process $\xi^{u}_{t}$ on
$G$ survives for a significant amount of time (call this event
\emph{quasi-survival}) should also be about $p_{\lambda}$.

The contact process is self-dual, in particular, the Poisson processes
used in the standard graphical construction can be reversed without
change of distribution. Thus, the event that $v \in \xi^{u}_{t}$ has
the same $\zz{P}_G$-probability as the event that $u\in \xi^{v}_{t}$, and
these events have the same $\zz{P}_G$-probability that two
\emph{independent} contact processes $\xi^{v}_{s}$ and $\xi^{u}_{s}$
started at $u$ and $v$ will intersect at time $t/2$. But
this will only happen if both  contact processes survive for time
$t/2$, and for large $t$ the probability of this will be about
$p_{\lambda}^{2}$. Hence, for large $t$, with high probability,
\begin{equation}\label{eq:upper}
	\zz{P}_G\{v\in \xi^{u}_{t} \}=\zz{P}_G\{ \xi^u_{t/2}\cap
	\xi^v_{t/2}\neq \emptyset\} \leq p_{\lambda}^{2} (1+o (1)).
\end{equation}

This argument shows that $p_{\lambda}^{2}$ is the largest possible
asymptotic value for the probability in relation~\eqref{eq:mainResultUpper}.
To show that this value is actually attained, we will show that
conditional on the
event of simultaneous quasi-survival for two independent contact
processes $\xi^{u}_{s}$ and $\xi^{v}_{s}$, the random sets
$\xi^{u}_{t_{+}/2}$ and $\xi^{v}_{t_{+}/2}$ will almost certainly
overlap. For this, we will argue that on the event of quasi-survival,
the cardinality of $\xi^{u}_{t_{+}/2}$  will be at least
$n^{1/2+\delta}$ for some $\delta >0$ depending on $\varepsilon$, and
that  $\xi^{v}_{t_{+}/2}$ is approximately distributed as a random
subset of $[n]$ of cardinality $n^{1/2-\delta}$. Since two such
independent random subsets will intersect with high probability, this
suggests that
\[
	\zz{P}_G\{ \xi^u_{t/2}\cap
	\xi^v_{t/2}\neq \emptyset\} \approx p_{\lambda}^{2}.
\]
Unfortunately, because the labeling processes used in constructing the
two contact processes $\xi^{v}_{s}$ and $\xi^{u}_{s}$  will interfere,
it will turn out that the  random sets $\xi^{u}_{t_{+}/2}$ and
$\xi^{v}_{t_{+}/2}$ are not independent, and so a more circuitous
argument will be needed.

Henceforth, let $t_1=(1-\varepsilon)\log n/2c_\lambda$,
$t_2=(1+3\varepsilon)\log n/2c_\lambda$ and $\Delta\,
t=t_2-t_1=2\varepsilon\log n/c_\lambda$. Notice that
$t_1+t_2=t_+$.  We will run $\xi^u$ for $t_2$ and
$\xi^v$ for $t_1$. In order to show
$\mathbb{P}\{\xi^u_{t_2}\cap \xi^v_{t_1}\neq \emptyset\}\geq
(1-o(1))p^{2}_{\lambda}$, we will apply the second moment
method. Let us briefly describe the idea.

Fix $0<\varepsilon<1/8$ and let $0<\delta\ll \varepsilon$. We will first grow
$\tilde{\xi}^u$ up to time $t_1$, and $\tilde{\xi}^v$ up to time $t_1-1$.
We call the processes to be in \emph{stage 1}. Conditional on
the event that both contact processes survive (which has probability
about $p^{\,2}_{\lambda}$), with high probability, the sets
of pioneer points,
$|\tilde{\zeta}^u_{t_1}|$ and $|\tilde{\zeta}^v_{t_1-1}|$ will be moderately
large (at least $n^{(1-\varepsilon)(1-\delta)/2}$), while $|\cup_{s\leq
t_1}\tilde{\xi}^u_s|$ and $|\cup_{s\leq t_1}\tilde{\xi}^v_s|$ will not be too
large (no more than $n^{(1-\varepsilon)(1+\delta)/2}$). Based on these,
with high probability, in the exploration process, we observe (1) the
labels assigned to $\cup_{s\leq t_1}\tilde{\xi}^u_s$ have no overlap
with the labels assigned to $\cup_{s\leq t_1-1}\tilde{\xi}^v_s$, and
(2) the labeling in $\cup_{s\leq t_1}\tilde{\xi}^u_s$ does not lead
to coexisting infections on the same site in the finite graph, and
thus $|\tilde{\xi}^u_s|=|\xi^u_s|$, for all $s\leq t_1$, and the same
holds for $|\tilde{\xi}^v_s|$ up to time $t_1-1$.
If all above events occur (which has probability $(1-o(1))p^{2}_{\lambda}$),
we call it a \emph{good event}.

Our goal is to show that conditional on the good event happening
in stage 1, if we run $\xi^u_s$ for another time $\Delta t$, and
run $\xi^v_s$ for another time $1$, then with high probability,
there will be at least one common label between them at the end.
We call this period to be \emph{stage 2}.
It suffices to consider certain subsets of both contact processes
which are much easier to deal with.
At the end of stage 1, we only keep the pioneer points of $\tilde{\xi}^u_{t_1}$.
For each such pioneer point (call it $i$), together with its free branch,
we run an independent severed contact process $(\tilde{\eta}^i_{s})_{s\geq 0}$
inside this branch for another duration of $\Delta t$, and label all of them.
We run a similar process for each pioneer point of $\tilde{\xi}^v_{t_1-1}$
(call it $j$, and the severed contact process $\tilde{\eta}^j_{\Delta t}$)
for time $1$, and label all of them. Such constructed $\phi(\cup_i \tilde{\eta}^i_{\Delta t})$
should be regarded as a subset of $\xi^u_{t_2}$, and  $\phi(\cup_j \tilde{\eta}^j_{1})$
should be regarded as a subset of $\xi^v_{t_1}$.

For each pair of $(i,j)$, let $I_{i,j}$
be the event that at the end of stage 2 we observe a common label assigned
between $\tilde{\eta}^i_{\Delta t}$ and $\tilde{\eta}^j_{1}$ in
a very specific way which we will state later. We will use a second moment
argument to show that $\sum_{i}\sum_{j}I_{ij}\rightarrow \infty$ with
high probability, which implies that $\xi^u_{t_2}\cap \xi^v_{t_1}\neq \emptyset$
with high probability.

Let us formulate the terminologies above in more details. We define $\goodevent$
(the good event mentioned above) to be the event that all of the following
events happen in stage 1:
\begin{compactenum}
\item [(1)] $|\tilde{\zeta}^u_{t_1}|\geq n^{(1-\varepsilon)(1-\delta)/2}$,
$|\tilde{\zeta}^v_{t_1-1}|\geq n^{(1-\varepsilon)(1-\delta)/2}$.

\item [(2)]$|\cup_{s\leq t_1}\tilde{\xi}^u_{s}|\leq n^{(1-\varepsilon)(1+\delta)/2}$,
$|\cup_{s\leq t_1-1}\tilde{\xi}^v_{s}|\leq n^{(1-\varepsilon)(1+\delta)/2}$.

\item [(3)]
In the exploration process of $|\cup_{s\leq t_1}\tilde{\xi}^u_{s}|$,
there is not a pair of distinct vertices in $|\cup_{s\leq t_1}\tilde{\xi}^u_{s}|$
on the cover tree that are assigned the same label on the finite graph.
The same holds for $|\cup_{s\leq t_1-1}\tilde{\xi}^v_{s}|$. Moreover,
the labels assigned to $|\cup_{s\leq t_1}\tilde{\xi}^u_{s}|$ are
completely different from these to $|\cup_{s\leq t_1-1}\tilde{\xi}^v_{s}|$.
\end{compactenum}

Notice that on $\goodevent$, all infections in $\tilde{\xi}^u_{t_1}$ and
$\tilde{\xi}^v_{t_1-1}$ are colored BLUE. Without loss of generality we will
label $\cup_{s\leq t_1}\tilde{\xi}^u_{s}$ first, then
$\cup_{s\leq t_1-1}\tilde{\xi}^v_{s}$.

\begin{Proposition}\label{Moment.Goodevent}
$\mathbb{P}\{\goodevent \}\geq(1-o(1))p^{2}_{\lambda}$, \text{\rm as }$n\rightarrow \infty$.
\end{Proposition}
\begin{proof}
For the two contact processes on the cover tree, with probability at
least $p^{2}_{\lambda}$, one survive up to time $t_1$ and the other up to
time $t_1-1$. Then it follows from Proposition \ref{CP.Uniongrowthrate}
and \ref{CP.pioneerpoint} that (1) and (2) in the definition of $\goodevent$
hold simultaneously with probability at least $(1-o(1))p^{2}_{\lambda}$.

Given (1) and (2), (3) is a statement on the exploration process,
and we will use the configuration model to estimate its probability.
When we first label the vertices in $\cup_{s\leq t_1}\tilde{\xi}^u_{s}$,
(3) requires that whenever there is a vertex to be labeled,
it cannot use any of those labels already used. For example, $u$ is
automatically the first label; for the second vertex, (3) requires it
not to be labeled $u$, and this will have probability $1-(2d-1)/(dn-1)$,
because an half edge with label $u$ is already used as the first half-edge,
and among the remaining $dn-1$ half edges, $d-1$ with label $u$ and $d$
with label $v$ are excluded by (3). Similarly, when we label the $(m+1)$-st
vertex, given the previous labeling is consistent with (3), the chance it
does not coincide with any of the previously used labels is $1-((d-2)m+d+1)/(dn-2m+1)$.

Therefore, the chance that the labeling in $|\cup_{s\leq t_1}\tilde{\xi}^u_{s}|$
being consistent with (3) is
\begin{equation}\label{Moment.Goodeventprob1}
\prod_{m=1}^{A-1} \left(1-\frac{(d-2)m+d+1}{dn-2m+1}\right),
\end{equation}
where $A=|\cup_{s\leq t_1}\tilde{\xi}^u_{s}|\leq n^{(1-\varepsilon)(1+\delta)/2}$.

Similarly, given the labeling in $|\cup_{s\leq t_1}\tilde{\xi}^u_{s}|$ is
consistent with (3), the chance that the labeling in
$|\cup_{s\leq t_1-1}\tilde{\xi}^v_{s}|$ comply with (3), is given by
\begin{equation}\label{Moment.Goodeventprob2}
\prod_{m=1}^{B-1} \left(1-\frac{(d-2)(A-1)+(d-2)(m-1)+1)}{dn-(2A-2)-2m+1}\right),
\end{equation}
where $B=|\cup_{s\leq t_1-1}\tilde{\xi}^v_{s}|\leq n^{(1-\varepsilon)(1+\delta)/2}$.

So the chance that (3) holds given (1) and (2) is the product of
(\ref{Moment.Goodeventprob1}) and (\ref{Moment.Goodeventprob2}),
which is no less than
\begin{equation*}\label{Moment.Goodeventprob}
\begin{aligned}
&\left(1-\frac{(d-2)(A+B-2)+1}{dn-o(n)}\right)^{A+B}\\
=&\exp\left(-O\left(\frac{(A+B)^2}{n}\right)\right)
\geq\exp\left(-O(n^{-\varepsilon-\delta+\varepsilon\delta})\right)\\
=& 1-o(1),\text{\rm as } n\rightarrow\infty,
\end{aligned}
\end{equation*}
and therefore the $\mathbb{P}$-probability of observing $\goodevent$
is at least $(1-o(1))p^{2}_{\lambda}$.
\end{proof}

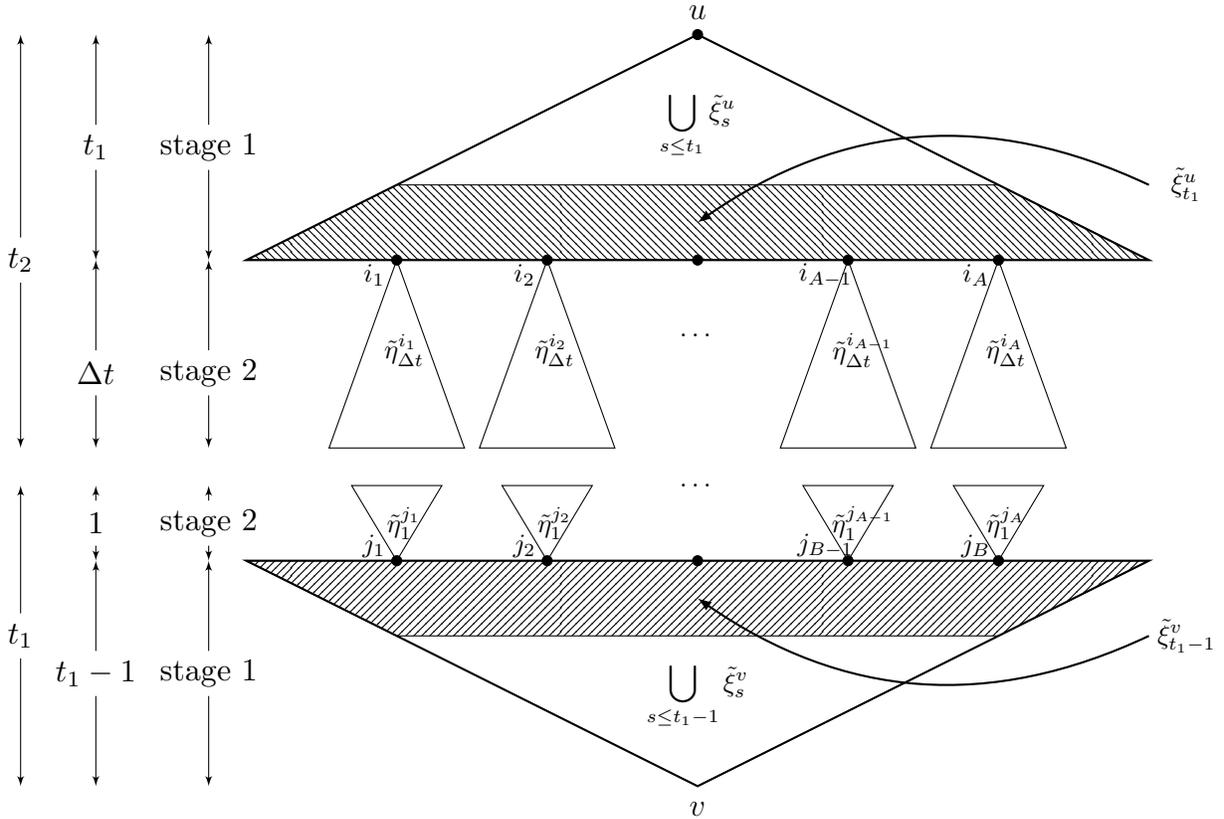
\begin{figure}[H]
\centering
\begin{tikzpicture}[scale=1]
\node[draw, circle, scale=0.4, fill=black] (u) at (0,0) {};
\node[] at (0,0.3) {\large $u$};
\draw[thick](0,0) -- (-6,-3) -- (6,-3) -- cycle;
\draw[pattern=north west lines] (-4,-2) -- (4,-2) -- (6, -3) -- (-6, -3) -- cycle;
\node[] at (0, -1.2) {$\displaystyle \bigcup_{s\le t_{1}} \tilde{\xi}_{s}^{u}$};
\draw[thick][->, bend right, -latex] (6,-2) to (0, -2.5);
\node[] at (6.5, -2) {$\displaystyle  \tilde{\xi}_{t_1}^{u}$};
\node[draw, circle, scale=0.4, fill=black] (i) at (-4, -3) {};
\node[]  at (-4.3, -3.2) {$i_1$};
\node[draw, circle, scale=0.4, fill=black] (i) at (-2, -3) {};
\node[]  at (-2.3, -3.2) {$i_2$};
\node[draw, circle, scale=0.4, fill=black] (i) at (2, -3) {};
\node[]  at (1.7, -3.2) {$i_{A-1}$};
\node[draw, circle, scale=0.4, fill=black] (i) at (4, -3) {};
\node[]  at (3.7, -3.2) {$i_{A}$};
\node[draw, circle, scale=0.4, fill=black] (i) at (0, -3) {};
\draw(-4,-3) -- (-4.9,-5.5) -- (-3.1,-5.5) -- cycle;
\node[scale=1] at (-3.9, -4.2) {$\displaystyle   \tilde{\eta}_{\Delta t}^{i_1}$};
\draw(-2,-3) -- (-2.9,-5.5) -- (-1.1,-5.5) -- cycle;
\node[scale=1] at (-1.9, -4.2) {$\displaystyle    \tilde{\eta}_{\Delta t}^{i_2}$};
\draw(2,-3) -- (2.9,-5.5) -- (1.1,-5.5) -- cycle;
\node[scale=1] at (2.2, -4.2) {$\displaystyle   \tilde{\eta}_{\Delta t}^{i_{A-1}}$};
\draw(4,-3) -- (4.9,-5.5) -- (3.1,-5.5) -- cycle;
\node[scale=1] at (4.1, -4.2) {$\displaystyle    \tilde{\eta}_{\Delta t}^{i_A}$};
\node[] at (0, -4) {$\dots$};
\node[draw, circle, scale=0.4, fill=black] (i) at (-4, -7) {};
\node[]  at (-4.3, -6.8) {$j_1$};
\node[draw, circle, scale=0.4, fill=black] (i) at (-2, -7) {};
\node[]  at (-2.3, -6.8) {$j_2$};
\node[draw, circle, scale=0.4, fill=black] (i) at (2, -7) {};
\node[]  at (1.7, -6.8) {$j_{B-1}$};
\node[draw, circle, scale=0.4, fill=black] (i) at (4, -7) {};
\node[]  at (3.7, -6.8) {$j_{B}$};
\node[draw, circle, scale=0.4, fill=black] (i) at (0, -7) {};
\draw[thick][->, bend left, -latex] (6,-8) to (0, -7.5);
\node[] at (6.5, -8) {$\displaystyle  \tilde{\xi}_{t_1-1}^{v}$};
\draw(-4,-7) -- (-4.6,-6) -- (-3.4,-6) -- cycle;
\node[scale=1] at (-3.9, -6.5) {$\displaystyle  \tilde{\eta}_{1}^{j_1}$};
\draw(-2,-7) -- (-2.6,-6) -- (-1.4,-6) -- cycle;
\node[scale=1] at (-1.9, -6.5) {$\displaystyle   \tilde{\eta}_{1}^{j_2}$};
\draw(2,-7) -- (2.6,-6) -- (1.4,-6) -- cycle;
\node[scale=1] at (2.2, -6.5) {$\displaystyle \tilde{\eta}_{1}^{j_{A-1}}$};
\draw(4,-7) -- (4.6,-6) -- (3.4,-6) -- cycle;
\node[scale=1] at (4.1, -6.5) {$\displaystyle   \tilde{\eta}_{1}^{j_A}$};
\node[] at (0, -10.3) {\large $v$};
\draw[thick] (0,-10) -- (-6,-7) -- (6,-7) -- cycle;
\draw[pattern=north east lines] (-6,-7) -- (6,-7) -- (4, -8) -- (-4, -8) -- cycle;
\node[] at (0, -8.8) {$\displaystyle \bigcup_{s\le t_{1}-1} \tilde{\xi}_{s}^{v}$};
\node[] at (0, -6) {$\dots$};
\node[scale=1.2] at (-8,-1.5) {$\displaystyle t_1$};
\draw[arrows=<->, -latex'] (-8,-1.8) -- (-8,-3);
\draw[arrows=<->, -latex'] (-8,-1.2) -- (-8,0);
\node[scale=1.2] at (-8,-4.5) {$\displaystyle \Delta t$};
\draw[arrows=<->, -latex'] (-8,-4.8) -- (-8,-5.5);
\draw[arrows=<->, -latex'] (-8,-4.2) -- (-8,-3);
\node[scale=1.2] at (-8,-8.5) {$\displaystyle t_1-1$};
\draw[arrows=<->, -latex'] (-8,-8.8) -- (-8,-10);
\draw[arrows=<->, -latex'] (-8,-8.2) -- (-8,-7);
\node[scale=1.2] at (-8,-6.5) {$\displaystyle 1$};
\draw[arrows=<->, -latex'] (-8,-6.8) -- (-8,-7);
\draw[arrows=<->, -latex'] (-8,-6.2) -- (-8,-6);
\node[scale=1.2] at (-6.5,-1.5) {stage 1};
\draw[arrows=<->, -latex'] (-6.5,-1.8) -- (-6.5,-3);
\draw[arrows=<->, -latex'] (-6.5,-1.2) -- (-6.5,0);
\node[scale=1.2] at (-6.5,-8.5) {stage 1};
\draw[arrows=<->, -latex'] (-6.5,-8.8) -- (-6.5,-10);
\draw[arrows=<->, -latex'] (-6.5,-8.2) -- (-6.5,-7);
\node[scale=1.2] at (-6.5,-4.5) {stage 2};
\draw[arrows=<->, -latex'] (-6.5,-4.8) -- (-6.5,-5.5);
\draw[arrows=<->, -latex'] (-6.5,-4.2) -- (-6.5,-3);
\node[scale=1.2] at (-6.5,-6.5) {stage 2};
\draw[arrows=<->, -latex'] (-6.5,-6.8) -- (-6.5,-7);
\draw[arrows=<->, -latex'] (-6.5,-6.2) -- (-6.5,-6);
\node[scale=1.2] at (-9,-3) {$\displaystyle t_2$};
\draw[arrows=<->, -latex'] (-9,-2.7) -- (-9,0);
\draw[arrows=<->, -latex'] (-9,-3.3) -- (-9,-5.5);
\node[scale=1.2] at (-9,-8) {$\displaystyle t_1$};
\draw[arrows=<->, -latex'] (-9,-7.7) -- (-9,-6);
\draw[arrows=<->, -latex'] (-9,-8.3) -- (-9,-10);
\end{tikzpicture}
\caption {$i_1,\dots,i_A$ are all pioneer points of $\tilde{\xi}_{t_1}^u$;
$j_1,\dots,j_B$ are all pioneer points of $\tilde{\xi}_{t_1-1}^v$.
We will run independent severed contact processes inside these
branches associated with the pioneer points.    \label{fig:t1} }
\end{figure}

From now on, we assume that in stage 1, $\goodevent$ happens.
Notice on $\goodevent$, there is a one-to-one correspondence between
infections on the cover tree and infections on the finite graph up to
time $t_1$ ($t_1-1$) for the contact processes started from $u$ ($v$).
For each pioneer point $i\in\zeta^u_{t_1}$ together with its free branch,
let $(\tilde{\eta}^i_{s})_{s\geq 0}$
denote the severed contact process in this branch with $\{ i \}$ being the initial
configuration. Similarly, for each pioneer point $j\in\zeta^v_{t_1-1}$ let
$(\tilde{\eta}^j_{s})_{s\geq 0}$ denote the corresponding severed contact process.
See Figure \ref{fig:t1} for an graphical illustration.

\begin{figure}[H]
\centering
\begin{tikzpicture}[scale=1]
\node[draw, circle, scale=0.4, fill=black] (u) at (0,0) {};
\node[draw, circle, scale=0.4, fill=black] (u) at (0,-9) {};
\node[] at (0,0.3) {\large $u$};
\node[] at (0.3,-2.2) {\large $i$};
\node[] at (0.2,-6.8) {\large $j$};
\draw[thick](0,0) -- (-6,-2) -- (6,-2) -- cycle;
\node[] at (0, -1.2) {$\displaystyle \bigcup_{s\le t_{1}} \tilde{\xi}_{s}^{u}$};
\node[draw, circle, scale=0.4, fill=black] (i) at (0, -2) {};
\node[draw, circle, scale=0.4, fill=black] (i) at (0, -7) {};
\node[] at (0, -9.3) {\large $v$};
\draw[thick](0,-2) -- (-1.5,-5) -- (1.5,-5) -- cycle;
\draw[pattern=north east lines] (-1.25,-4.5) -- (1.25,-4.5) -- (1.5, -5) -- (-1.5, -5) -- cycle;
\draw[thick][->, bend right, -latex] (3,-5) to (0, -4.75);
\node[] at (3.5,-5) {\large $\tilde{\eta}^i_{\Delta t -1}$};
\node[draw, circle, scale=0.4, fill=black] at (0,-5) {};
\node[] at (0, -3.8) {$\displaystyle \bigcup_{s\le \Delta t-1} \tilde{\eta}_{s}^{i}$};
\node[] at (0.3, -5.2) {\large $x_0$};
\draw[thick] (0,-9) -- (-6,-7) -- (6,-7) -- cycle;
\node[draw, circle, scale=0.4, fill=black] at (0,-6) {};
\node[] at (0.7, -6) {\large $x_1=y$};
\draw[thick] (0,-5) -- (0,-6) -- (0,-7) -- cycle;
\node[] at (0, -7.8) {$\displaystyle \bigcup_{s\le t_{1}-1} \tilde{\xi}_{s}^{v}$};
\node[scale=1.2] at (-8,-3.5) {$\displaystyle \Delta t-1$};
\draw[arrows=<->, -latex'] (-8,-3.7) -- (-8,-5);
\draw[arrows=<->, -latex'] (-8,-3.3) -- (-8,-2);
\node[scale=1.2] at (-8,-5.5) {$\displaystyle 1$};
\draw[arrows=<->, -latex'] (-8,-5.3) -- (-8,-5);
\draw[arrows=<->, -latex'] (-8,-5.7) -- (-8,-6);
\node[scale=1.2] at (-8,-6.5) {$\displaystyle 1$};
\draw[arrows=<->, -latex'] (-8,-6.7) -- (-8,-7);
\draw[arrows=<->, -latex'] (-8,-6.3) -- (-8,-6);
\node[scale=1.2] at (-6.5,-1) {stage 1};
\draw[arrows=<->, -latex'] (-6.5,-1.2) -- (-6.5,-2);
\draw[arrows=<->, -latex'] (-6.5,-0.8) -- (-6.5,0);
\node[scale=1.2] at (-6.5,-8) {stage 1};
\draw[arrows=<->, -latex'] (-6.5,-8.2) -- (-6.5,-9);
\draw[arrows=<->, -latex'] (-6.5,-7.8) -- (-6.5,-7);
\node[scale=1.2] at (-6.5,-4) {stage 2};
\draw[arrows=<->, -latex'] (-6.5,-4.2) -- (-6.5,-6);
\draw[arrows=<->, -latex'] (-6.5,-3.8) -- (-6.5,-2);
\node[scale=1.2] at (-6.5,-6.5) {stage 2};
\draw[arrows=<->, -latex'] (-6.5,-6.3) -- (-6.5,-6);
\draw[arrows=<->, -latex'] (-6.5,-6.7) -- (-6.5,-7);
\end{tikzpicture}
\caption{
The event $I_{ij}$. $x_0$ is a pioneer point of $\tilde{\eta}^i_{\Delta t -1}$.
\label{fig:xy}
}
\end{figure}
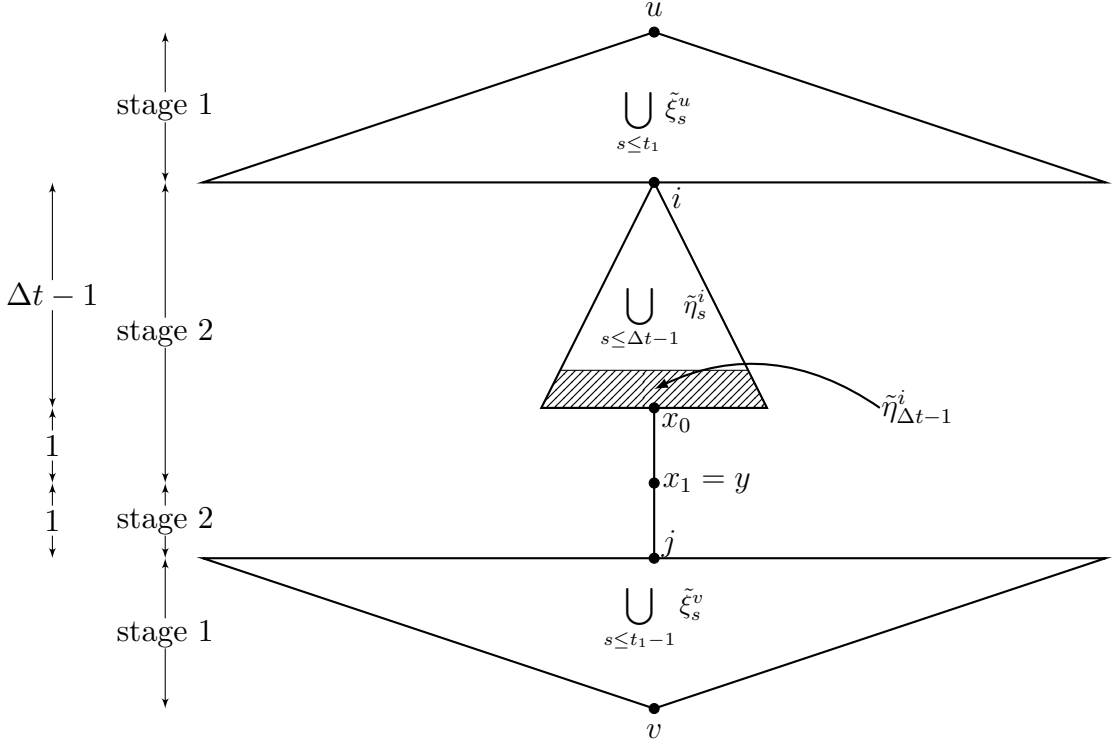

Next we will define the event $I_{ij}$ for $i\in \tilde{\zeta}^u_{t_1}$ and
$j\in \tilde{\zeta}^v_{t_1-1}$. We say $I_{i,j}$ happens if in stage 2
all of the following events happen:
\begin{compactenum}
\item[(1)] for $(\tilde{\eta}^j_{s})_{0\leq s \leq 1}$, $j$ infects
the neighbor in its free branch, call it $y$, before time $1$,
and the infection at $y$ stays alive but does not give rise to
an infection till the end of time $1$;
\item [(2)] for $(\tilde{\eta}^i_{s})_{0\leq s \leq \Delta t}$, up to
time $\Delta t $, labels assigned to vertices in
$\cup_{s\leq \Delta t}\tilde{\eta}^i_{s}$ do not appear in the
label set assigned to $\cup_{s\leq t_1}\tilde{\xi}^u_s$,
$\cup_{s\leq t_1-1}\tilde{\xi}^v_s$, or any $\cup_{s\leq \Delta t}\tilde{\eta}^{i\prime}_s$
for $i^\prime \in \tilde{\zeta}^u_{t_1}, i^\prime \neq i$;
also, distinct vertices in $\cup_{s\leq \Delta t}\tilde{\eta}^i_{s}$
are assigned distinct labels;
\item [(3)]$x_0$, a pioneer point of $\tilde{\eta}^i_{\Delta t-1}$, infects
the neighbor in its free branch (call this neighbor $x_1$) in
time interval $[\Delta t-1,\Delta t]$,  and the
infection at $x_1$ stays alive till the end of $\Delta t$ without infecting other vertices;
furthermore, $x_1$ is assigned the same label as $y$.
\end{compactenum}
See Figure \ref{fig:xy} for a graphical illustration.

Now in order to show with probability approaching 1, there exists some
$i$ and $j$ such that $I_{i,j}$ happens, we will show the following
assertion on the first moment
\begin{equation} \tag{FM}\label{Moment.1stmoment}
\begin{aligned}
 &\sum_{i\in \tilde{\zeta}^u_{t_1}}\sum_{j\in \tilde{\zeta}^v_{t_1-1}}\mathbb{P}\{I_{ij}| \goodevent \}  \rightarrow \infty, \text{\rm as } n\rightarrow\infty , \\
 \end{aligned}
\end{equation}
and the following assertion on the second moment
\begin{equation}\tag{SM}\label{Moment.2ndmoment}
\begin{aligned}
&\sum_{i\in \tilde{\zeta}^u_{t_1}}\sum_{j\in \tilde{\zeta}^v_{t_1-1}}
\sum_{i^\prime \in \tilde{\zeta}^u_{t_1}}\sum_{j^\prime \in \tilde{\zeta}^v_{t_1-1}}
\mathbb{P}\{I_{ij}\cap I_{i^\prime j^\prime}| \goodevent \} =(1+o(1) )  \left(\sum_{i\in \tilde{\zeta}^u_{t_1}}\sum_{j\in \tilde{\zeta}^v_{t_1-1}}\mathbb{P}\{I_{ij}| \goodevent\}  \right)^2,  \text{\rm as } n\rightarrow\infty .
\end{aligned}
\end{equation}
The next two subsections are devoted to proving (\ref{Moment.1stmoment}) and (\ref{Moment.2ndmoment}).

\subsection{1st Moment Calculation}

\begin{proof}[Proof of (\ref{Moment.1stmoment})]
Since $\mathbb{P}\{I_{ij}|\goodevent\}$ is constant among all pairs
$(i,j)\in \tilde{\zeta}^u_{t_1}\times \tilde{\zeta}^v_{t_1-1}$, it suffices to
estimate a single term.
Let us first estimate the size of $|\cup_{s\leq \Delta t}\tilde{\eta}^i_s|$.
From Proposition \ref{CP.exponentialbound},
\begin{equation*}\label{2ndmoment.10}
\zz{P}\{|\cup_{s\leq \Delta t}\tilde{\eta}^i_s|>n^{2(1+\delta)\varepsilon}\}\leq \exp(-K \gamma^{\Delta t} )=\exp(-K \gamma^{2\varepsilon \log n/c_\lambda} ),
\end{equation*}
which is decaying faster than any polynomial of $1/n$.
On the other hand, on $\goodevent$, there are no more than
$n^{(1-\varepsilon)(1+\delta)/2}$ different $i$'s. The same
argument works for $j$'s. Therefore by a union bound,
with probability $1-o(1)$, we have
\begin{equation}\label{Moment.uniformbound}
|\cup_{s\leq \Delta t}\tilde{\eta}^i_s|\leq n^{2(1+\delta)\varepsilon},
\text{\rm for all }i\in\tilde{\zeta}^u_{t_1};\,
|\cup_{s\leq 1}\tilde{\eta}^j_s|\leq n^{2(1+\delta)\varepsilon},
\text{\rm for all }j\in\tilde{\zeta}^v_{t_1-1}.
\end{equation}

Let us estimate how many half-edges will be used at the
end of stage 2. If (\ref{Moment.uniformbound}) holds,
\begin{compactenum}
\item[(1)] labeling $\cup_{s\leq t_1}\tilde{\xi}^u_s$ and $\cup_{s\leq t_1-1}\tilde{\xi}^v_s$
have used no more than $O(n^{(1-\varepsilon)(1+\delta)/2})$ half-edges;
\item[(2)] labeling $\cup_{s\leq \Delta t}\tilde{\eta}^{i^\prime}_s$ and
$\cup_{s\leq \Delta t}\tilde{\eta}^{j^\prime}_s$ for all
$i^\prime\neq i$ and $j^\prime \neq j$ will use
no more than $O(n^{(1-\varepsilon)(1+\delta)/2}\times n^{2(1+\delta)\varepsilon})$,
which is $O(n^{(1+\delta)(1+3\varepsilon)/2})$ half-edges;
\item[(3)] labeling $\cup_{s\leq \Delta t}\tilde{\eta}^{i}_s$ and
$\cup_{s\leq 1}\tilde{\eta}^{j}_s$ will use no more than
$O(n^{2(1+\delta)\varepsilon})$ half-edges.
\end{compactenum}
Therefore we see that during the labeling procedure, at any time
there are at least $dn-O(n^{(1+3\varepsilon)(1+\delta)/2})$ unused half-edges
in the pool. Therefore, from the above analysis, a counting argument
similar to the proof of Proposition \ref{Moment.Goodevent} will conclude that
with probability $1-o(1)$, the labels assigned to
$\cup_{s\leq \Delta t}\tilde{\eta}^{i}_s$ and
$\cup_{s\leq 1}\tilde{\eta}^{j}_s$ do not overlap with
those assigned to $\cup_{s\leq t_1}\tilde{\xi}^u_s$,
$\cup_{s\leq t_1-1}\tilde{\xi}^v_s$,
$\cup_{s\leq \Delta t}\tilde{\eta}^{i^\prime}_s$ and
$\cup_{s\leq 1}\tilde{\eta}^{j^\prime}_s$ for all
$i^\prime\neq i$ and $j^\prime \neq j$; moreover, different
vertices in $\cup_{s\leq \Delta t}\tilde{\eta}^{i}_s$ are assigned
different labels; the same holds for $\cup_{s\leq 1}\tilde{\eta}^{j}_s$.

Now let us discuss how to produce a common label in the prescribed
way in the definition of $I_{ij}$.
Let $p_{\lambda,\text{severed}}$ be the chance that the severed
contact process with rate $\lambda$ survives. It is not hard to
deduce that $p_{\lambda,\text{severed}}>0$ using proposition
\ref{CP.severedgrowthrate}, and thus
\[\mathbb{P}\{\tilde{\eta}^i_{\Delta t-1}\neq \emptyset\}\geq p_{\lambda,\text{severed}}>0.\]
From Proposition \ref{CP.severedpioneerpoint}, with probability
$(1-o(1))p_{\lambda,\text{severed}}$,
$|\tilde{\psi}^i_{\Delta t-1}|\geq e^{-c_\lambda} n^{2(1-\delta)\varepsilon}$.

With some positive probability $q_\lambda>0$, for
$(\tilde{\eta}^j_{s})_{0\leq s\leq 1}$, $j$ infects the neighbor
in its free branch (call this neighbor $y$) before time $1$,
and the infection at $y$ stays alive but does not give rise to
an infection till the end of time $1$. Fix a vertex $x_0 \in
\tilde{\eta}^i_{\Delta t-1}$ which is a pioneer point, then
$q_\lambda$ is also the probability that the infection at $x_0$
infects the neighbor in its free branch (call this neighbor $x_1$)
before time $1$, and the infection at $x_1$ stays alive but does
not give rise to an infection till the end of time $1$. See figure \ref{fig:xy}.

As long as when we label $y$, we choose a label which has been
chosen for the first time (which is of probability $1-o(1)$, because the
number of used half-edges is of order $O(n^{(1+\delta)(1+3\varepsilon)/2})$),
then when we label $x_1$, the chance of choosing the same label as
$y$ is at least $(d-1)/dn$. Now we can apply a Binomial type of argument
to show that conditional on $|\tilde{\psi}^i_{\Delta t-1}|\geq e^{-c_\lambda} n^{2(1-\delta)\varepsilon}$, the chance of observing a common label
in the prescribed manner is at least
\[\frac{1}{2} (1-o(1)) e^{-c_\lambda} n^{2(1-\delta)\varepsilon}\frac{d-1}{dn}q_\lambda^2=C_1 n^{2(1-\delta)\varepsilon -1},\]
for some $C_1>0$.

Combining all the calculation above, we conclude that
\[\mathbb{P}\{I_{ij}|F\} \geq (1-o(1)) p_{\lambda,\text{severed}}C_1n^{2(1-\delta)\varepsilon -1}\geq C_2n^{2(1-\delta)\varepsilon -1},
\]
for some $C_2>0$.
Since there are at least $n^{(1-\varepsilon)(1-\delta)/2}$ distinct $i$'s
and $n^{(1-\varepsilon)(1-\delta)/2}$ distinct $j$'s on $\goodevent$,
we conclude
\[\sum_{i\in \tilde{\zeta}^u_{t_1}} \sum_{j\in \tilde{\zeta}^v_{t_1-1}}
\mathbb{P}\{I_{ij}|\goodevent\}\geq n^{(1-\varepsilon)(1-\delta)/2}
n^{(1-\varepsilon)(1-\delta)/2} C_2n^{2(1-\delta)\varepsilon -1}=C_2n^{\varepsilon-\delta(1+\varepsilon)}\rightarrow
\infty, \text{ as }n\rightarrow \infty,
\]
which is (\ref{Moment.1stmoment}).
\end{proof}

\subsection{2nd Moment Calculation}
This subsection is dedicated to showing (\ref{Moment.2ndmoment}).
We expand the second moment as follows,
\begin{equation*}\label{2ndmoment.expansion}
\begin{aligned}
&\sum_{i\in \tilde{\zeta}^u_{t_1}}\sum_{j\in \tilde{\zeta}^v_{t_1-1}}
\sum_{i^\prime \in \tilde{\zeta}^u_{t_1}}\sum_{j^\prime \in \tilde{\zeta}^v_{t_1-1}}
\mathbb{P}\{I_{ij}\cap I_{i^\prime j^\prime}| \goodevent \}
=I+II+III+IV,\\
\end{aligned}
\end{equation*}
where
\begin{equation*}
\begin{aligned}
I&=\sum_{i\in \tilde{\zeta}^u_{t_1}}\sum_{j\in \tilde{\zeta}^v_{t_1-1}}
\mathbb{P}\{I_{ij}| \goodevent \},  \\
II&=\sum_{i\in \tilde{\zeta}^u_{t_1}}\sum_{j\neq j^\prime\in \tilde{\zeta}^v_{t_1-1}}
\mathbb{P}\{I_{ij}\cap I_{ij^\prime}| \goodevent \}, \\
III&=\sum_{i\neq i^\prime \in \tilde{\zeta}^u_{t_1}}\sum_{j\in \tilde{\zeta}^v_{t_1-1}}
\mathbb{P}\{I_{ij}\cap I_{i^\prime j}| \goodevent \}, \\
IV&=\sum_{i\neq i^\prime \in \tilde{\zeta}^u_{t_1}}\sum_{j\neq j^\prime \in \tilde{\zeta}^v_{t_1-1}}
\mathbb{P}\{I_{ij}\cap I_{i^\prime j^\prime}| \goodevent \}.
\end{aligned}
\end{equation*}
Now (\ref{Moment.2ndmoment}) becomes
\[\label{2ndmoment.simplification}
I+II+III+IV=\left(1+o(1)\right)I^2, \text{\rm as }n\rightarrow \infty.
\]
We have already shown that $I\rightarrow \infty$ as $n\rightarrow \infty$, so
\[\label{2ndmoment.I finished}
I=o(I^2).
\]
We will show the following assertions in the next subsections,
\begin{align}
II&=o(I^2), \label{2ndmoment.II finished} \\
III&=o(I^2), \label{2ndmoment.III finished} \\
IV&=(1+o(1))I^2. \label{2ndmoment.IV finished}
\end{align}

\subsubsection{Proof of (\ref{2ndmoment.II finished}) and (\ref{2ndmoment.III finished})}
\begin{proof}
We will only prove (\ref{2ndmoment.II finished}) because the proof of
(\ref{2ndmoment.III finished}) is essentially the same.
We will find an upper bound for
\[
\mathbb{P}\{I_{ij}\cap I_{ij^\prime}| \goodevent \}.
\]
From the proof of (\ref{Moment.uniformbound}) we know that
with probability at least  $1-n^{-4}$, (\ref{Moment.uniformbound}) holds.
We condition on (\ref{Moment.uniformbound}) happening.
Suppose $y$ is the neighbor of $j$ in $j$'s free branch, and $y^\prime$
is the neighbor of $j^\prime$ in $j^\prime$'s free branch. In order that
$I_{ij}\cap I_{i j^\prime}$ happens, at least we should observe that when
we label $y$ we assign some label that has already been used in
$\cup_{s\leq \Delta t}\tilde{\eta}^i_s$, and the same is true for $y^\prime$.
$y$ has chance no more than
\[\frac{(d-1)n^{2(1+\delta)\varepsilon}}{dn-O(n^{(1+\delta)(1+3\varepsilon)/2})}\]
of using a label already used in $\cup_{s\leq \Delta t}\tilde{\eta}^i_s$,
because the numerator is an upperbound of qualified half-edges in the
pool, and the denominator is an lower bound of the number of unused
half-edges at the end of stage 2 if (\ref{Moment.uniformbound}) holds.
The same holds for $y^\prime$.

Therefore we obtain that
\[\mathbb{P}\{I_{ij}\cap I_{i j^\prime}\}\leq \frac{1}{n^4}+\left(\frac{(d-1)n^{2(1+\delta)\varepsilon}}{dn-O(n^{(1+\delta)(1+3\varepsilon)/2})}\right)^2=O(n^{4(1+\delta)\varepsilon -2}).\]

On $\goodevent$, there are no more than $n^{(1-\varepsilon)(1+\delta)/2}$
different $i$'s, and no more than $n^{(1-\varepsilon)(1+\delta)/2}$ different
$j$'s and hence no more than $n^{(1-\varepsilon)(1+\delta)}$ different pairs
of $(j,j^\prime)$. Therefore we sum over all $i,j,j^\prime$ and obtain
\[
\begin{aligned}
\sum_{i\in \tilde{\zeta}^u_{t_1}}\sum_{j\neq j^\prime\in \tilde{\zeta}^v_{t_1-1}}
\mathbb{P}\{I_{ij}\cap I_{ij^\prime}| \goodevent \}
\leq n^{(1-\varepsilon)(1+\delta)/2} n^{(1-\varepsilon)(1+\delta)} O(n^{4(1+\delta)\varepsilon -2})
=o(I^2).
\end{aligned}
\]
Therefore we've proved (\ref{2ndmoment.II finished}).
\end{proof}

\subsubsection{Proof of (\ref{2ndmoment.IV finished})}

For $i\neq i^\prime \in \tilde{\zeta}^u_{t_1}$ and
$j\neq j^\prime \in \tilde{\zeta}^v_{t_1-1}$ we will estimate
\[\mathbb{P}\{I_{ij}\cap I_{i^\prime j^\prime}|\goodevent \}. \]
Our goal is to show the above quantity is almost the same
as
\[\mathbb{P}\{I_{ij}|\goodevent\} \mathbb{P}\{I_{i^\prime j^\prime}|\goodevent \}. \]
From the proof of (\ref{Moment.uniformbound}) we know that
with probability at least  $1-n^{-4}$, (\ref{Moment.uniformbound}) holds.

First of all, the severed contact processes on the cover trees, namely
$\tilde{\eta}^i_{\Delta t}$, $\tilde{\eta}^{i^\prime}_{\Delta t}$,
$\tilde{\eta}^j_{1}$ and $\tilde{\eta}^{j^\prime}_{1}$ are independent.
Given a realization of the above 4 processes, now we consider
the labeling process on them. Let $A$ be a specific pattern of labeling
$\cup_{s\leq \Delta t} \tilde{\eta}^i_s$ and $\cup_{s\leq 1} \tilde{\eta}^j_s$,
and let $B$ be a specific pattern of labeling $\cup_{s\leq \Delta t} \tilde{\eta}^{i^\prime}_s$
and $\cup_{s\leq 1} \tilde{\eta}^{j^\prime}_s$.

Notice that
\begin{equation} \label{2ndmoment.expression1}
\begin{aligned}
\mathbb{P}\{I_{ij}|\goodevent\}&=
\sum_{A \text{ compatible with }I_{ij}}\mathbb{P}\{\text{pattern }A|\goodevent, \cup_{s\leq \Delta t} \tilde{\eta}^i_s,\cup_{s\leq   1} \tilde{\eta}^{j}_s\}
\mathbb{P}\{\cup_{s\leq \Delta t} \tilde{\eta}^i_s,\cup_{s\leq   1} \tilde{\eta}^{j}_s\}, \\
\mathbb{P}\{I_{i^\prime j^\prime}|\goodevent\}&=
\sum_{B \text{ compatible with }I_{i^\prime j^\prime}}\mathbb{P}\{\text{pattern }B|\goodevent, \cup_{s\leq \Delta t} \tilde{\eta}^{i^\prime}_s,\cup_{s\leq   1} \tilde{\eta}^{j^\prime}_s\}
\mathbb{P}\{\cup_{s\leq \Delta t} \tilde{\eta}^{i^\prime}_s, \cup_{s\leq   1} \tilde{\eta}^{j^\prime}_s\}, \\
\end{aligned}
\end{equation}
where in expressions such as $\mathbb{P}\{\cup_{s\leq \Delta t}
\tilde{\eta}^i_s,\cup_{s\leq   1} \tilde{\eta}^{j}_s\}$,
we only care about information that is relevant to the labeling
process, such as which vertices are ever infected and
the relative order of appearances of infections. Irrelevant
information such as when exactly an infection appears
is not included in this probability.

From now on let us use pattern $A\& B$ to denote pattern $A$ on
$\cup_{s\leq \Delta t} \tilde{\eta}^i_s,\cup_{s\leq   1} \tilde{\eta}^{j}_s$
and pattern $B$ on $\cup_{s\leq \Delta t} \tilde{\eta}^{i^\prime}_s,\cup_{s\leq   1} \tilde{\eta}^{j^\prime}_s$.
Also notice that if pattern $A$ is compatible with $I_{ij}$ and pattern $B$
is compatible with $I_{i^\prime j^\prime}$ then the  pattern $A\& B$
is compatible with $I_{ij}\cap I_{i^\prime j^\prime}$.
Therefore we can also express
\begin{equation}\label{2ndmoment.expression2}
 \begin{aligned}
&\mathbb{P}\{I_{ij}\cap I_{i^\prime j^\prime}|\goodevent\}=
\sum_{A \text{ compatible with }I_{ij} \,}
\sum_{B \text{ compatible with } I_{i^\prime j^\prime }} \\
&
\mathbb{P}\{\text{pattern } A\& B
\ \vert \goodevent, \cup_{s\leq \Delta t} \tilde{\eta}^i_s,\cup_{s\leq   1} \tilde{\eta}^{j}_s, \cup_{s\leq \Delta t} \tilde{\eta}^{i^\prime}_s,\cup_{s\leq   1} \tilde{\eta}^{j^\prime}_s\}
 \times \mathbb{P}\{\cup_{s\leq \Delta t} \tilde{\eta}^i_s,\cup_{s\leq   1} \tilde{\eta}^{j}_s)
\mathbb{P}\{\cup_{s\leq \Delta t} \tilde{\eta}^{i^\prime}_s,\cup_{s\leq   1} \tilde{\eta}^{j^\prime}_s\}.
\end{aligned}
\end{equation}
From the above expressions, it suffices to compare
\begin{equation} \label{2ndmoment.product}
\mathbb{P}\{\text{pattern }A|\goodevent, \cup_{s\leq \Delta t} \tilde{\eta}^i_s,\cup_{s\leq   1} \tilde{\eta}^{j}_s\} \times
\mathbb{P}\{\text{pattern }B|\goodevent, \cup_{s\leq \Delta t} \tilde{\eta}^{i^\prime}_s,\cup_{s\leq   1} \tilde{\eta}^{j^\prime}_s\}
\end{equation}
and
\begin{equation}\label{2ndmoment.joint}
\mathbb{P}\{\text{pattern }A \& B
\ \vert \goodevent, \cup_{s\leq \Delta t} \tilde{\eta}^i_s,\cup_{s\leq   1} \tilde{\eta}^{j}_s, \cup_{s\leq \Delta t} \tilde{\eta}^{i^\prime}_s,\cup_{s\leq   1} \tilde{\eta}^{j^\prime}_s\}.
\end{equation}
We can rewrite (\ref{2ndmoment.joint}) as
\begin{equation}\label{2ndmoment.joint2}
\begin{aligned}
&\mathbb{P}\{\text{pattern }A
\ \vert \goodevent, \cup_{s\leq \Delta t} \tilde{\eta}^i_s,\cup_{s\leq   1} \tilde{\eta}^{j}_s, \cup_{s\leq \Delta t} \tilde{\eta}^{i^\prime}_s,\cup_{s\leq   1} \tilde{\eta}^{j^\prime}_s\} \\
&\times
\mathbb{P}\{\text{pattern }B
\ \vert \goodevent, \cup_{s\leq \Delta t} \tilde{\eta}^i_s,\cup_{s\leq   1} \tilde{\eta}^{j}_s, \cup_{s\leq \Delta t} \tilde{\eta}^{i^\prime}_s,\cup_{s\leq   1} \tilde{\eta}^{j^\prime}_s,\text{ pattern }A\}.
\end{aligned}
\end{equation}
The following proposition compares (\ref{2ndmoment.product})
and  (\ref{2ndmoment.joint2}).
\begin{Proposition}\label{2ndmoment.comparison}
If (\ref{Moment.uniformbound}) happens, then
\begin{equation}\label{2ndmoment.ineq1}
\mathbb{P}\{\text{\rm pattern }A
\ \vert \goodevent, \cup_{s\leq \Delta t} \tilde{\eta}^i_s,\cup_{s\leq   1} \tilde{\eta}^{j}_s, \cup_{s\leq \Delta t} \tilde{\eta}^{i^\prime}_s,\cup_{s\leq   1} \tilde{\eta}^{j^\prime}_s\}
\leq (1+o(1)) \mathbb{P}\{\text{\rm pattern }A|\goodevent, \cup_{s\leq \Delta t} \tilde{\eta}^i_s,\cup_{s\leq   1} \tilde{\eta}^{j}_s\},
\end{equation}
and
\begin{equation}\label{2ndmoment.ineq2}
\begin{aligned}
&\mathbb{P}\{\text{\rm pattern }B
\ \vert \goodevent, \cup_{s\leq \Delta t} \tilde{\eta}^i_s,\cup_{s\leq   1} \tilde{\eta}^{j}_s, \cup_{s\leq \Delta t} \tilde{\eta}^{i^\prime}_s,\cup_{s\leq   1} \tilde{\eta}^{j^\prime}_s,\text{
\rm pattern }A\} \\
\leq &(1+o(1))\mathbb{P}\{\text{\rm pattern }B|\goodevent, \cup_{s\leq \Delta t} \tilde{\eta}^{i^\prime}_s,\cup_{s\leq   1} \tilde{\eta}^{j^\prime}_s\}  .
\end{aligned}
\end{equation}
\end{Proposition}
\begin{proof}
Let us only show (\ref{2ndmoment.ineq2}) since the proof of
(\ref{2ndmoment.ineq1}) is the same.
The difference between
\[\mathbb{P}\{\text{\rm pattern }B
\ \vert \goodevent, \cup_{s\leq \Delta t} \tilde{\eta}^i_s,\cup_{s\leq   1} \tilde{\eta}^{j}_s, \cup_{s\leq \Delta t} \tilde{\eta}^{i^\prime}_s,\cup_{s\leq   1} \tilde{\eta}^{j^\prime}_s,\text{
\rm pattern }A\}\]
and
\[
\mathbb{P}\{\text{\rm pattern }B|\goodevent, \cup_{s\leq \Delta t} \tilde{\eta}^{i^\prime}_s,\cup_{s\leq   1} \tilde{\eta}^{j^\prime}_s\}
\]
is that whenever we draw an unused half-edge from the pool, the number
of unused half-edge is reduced because some are already used in pattern
$A$. However if (\ref{Moment.uniformbound}) happens, then
\begin{compactenum}
\item[(1)] labeling pattern $A$ uses no more than
$O(n^{2(1+\delta)\varepsilon})$ half-edges;
\item[(2)] the total number of half-edges to be drawn in stage 2
is no more than $O(n^{(1+\delta)(1+3\varepsilon)/2})$;
\item[(3)] at any time in stage 2, the number of unused half-edges
is at least $K(n)$, where
\[K(n)=dn-O(n^{(1+\delta)(1+3\varepsilon)/2}).\]
\end{compactenum}
Therefore
\[\mathbb{P}\{\text{\rm pattern }B
\vert \goodevent, \cup_{s\leq \Delta t} \tilde{\eta}^i_s,\cup_{s\leq   1} \tilde{\eta}^{j}_s, \cup_{s\leq \Delta t} \tilde{\eta}^{i^\prime}_s,\cup_{s\leq   1} \tilde{\eta}^{j^\prime}_s,\text{
\rm pattern }A\}\]
is no more than
\[
\left(\frac{K(n)}{K(n)-O(n^{2(1+\delta)\varepsilon})}\right)^{O(n^{(1+\delta)(1+3\varepsilon)/2})}
\mathbb{P}\{\text{\rm pattern }B|\goodevent, \cup_{s\leq \Delta t} \tilde{\eta}^{i^\prime}_s,\cup_{s\leq   1} \tilde{\eta}^{j^\prime}_s\},
\]
which is
\[(1+o(1))\mathbb{P}\{\text{\rm pattern }B|\goodevent, \cup_{s\leq \Delta t} \tilde{\eta}^{i^\prime}_s,\cup_{s\leq   1} \tilde{\eta}^{j^\prime}_s\}.\]
\end{proof}

\begin{proof}[Proof of (\ref{2ndmoment.IV finished})]
Based on whether (\ref{Moment.uniformbound}) happen or not, by
(\ref{2ndmoment.expression1}), (\ref{2ndmoment.expression2}) and
Proposition \ref{2ndmoment.comparison},
\[
\mathbb{P}\{I_{ij}\cap I_{i^\prime j^\prime}|\goodevent \}\leq \frac{1}{n^4}+(1+o(1)\}
\mathbb{P}\{I_{ij}|\goodevent \} \mathbb{P}\{ I_{i^\prime j^\prime}|\goodevent \},
\]
therefore when we sum over all $i,i^\prime,j,j^\prime$,
we obtain
\[
\sum_{i\neq i^\prime \in \tilde{\zeta}^u_{t_1}}\sum_{j\neq j^\prime \in \tilde{\zeta}^v_{t_1-1}}
\mathbb{P}\{I_{ij}\cap I_{i^\prime j^\prime}| \goodevent\}\leq o(1)+(1+o(1))\left( \sum_{i \in \tilde{\zeta}^u_{t_1}}\sum_{j \in \tilde{\zeta}^v_{t_1-1}}
\mathbb{P}\{I_{ij}| \goodevent \}\right)^2,
\]
which is  (\ref{2ndmoment.IV finished}).
\end{proof}

\section{Asymptotic infection density}
Throughout this section we let $t_+ = (1 + \varepsilon) \log n/c_\lambda$.

Fix $0<\varepsilon<1/8$ and let $g_n(\varepsilon)$ be in Theorem \ref{THM1}.
We say that a pair of vertices $(u,v)\in [n]\times [n]$ is \emph{good} if
\[\zz{P}_G\{v\in \xi^u_{t_+}\}\geq(1-g_n(\varepsilon))p^2_{\lambda}.\]
We say that a vertex $u\in [n]$ is \emph{good} if the set $\{v\in
[n]: v\neq u, (u,v) \text{ \rm  is a good pair}\}$ has cardinality at
least $(1-\sqrt[4]{g_n(\varepsilon)})(n-1)$. Using Markov inequality it
is easy to deduce from Theorem \ref{THM1} that
\begin{Proposition}\label{2.typical}
For asymptotically almost every $G$, the number of good pairs is at least
$(1-\sqrt{g_n(\varepsilon)})n(n-1)$, and the number of good
vertices is at least $(1-\sqrt[4]{g_n(\varepsilon)})n$.
\end{Proposition}
The choice of $\sqrt{g_n}$ and $\sqrt[4]{g_n}$ in the definition
of good pair\slash vertex and in the above proposition is not crucial;
we only need them to be $o(1)$ terms. The above proposition
shows that such defined good pairs\slash vertices are indeed typical.

The next proposition states that a contact process started from a good vertex has decent chance to survive time $t_+$.

\begin{Proposition}\label{2.typical survival}
Suppose $u\in [n]$ is a good vertex, and $\xi^u_t$ is a contact
process with initial state $\{u\}$ on $G$. Then there exist constants
$h_n(\varepsilon)\rightarrow 0$ as $n\rightarrow \infty$ such
that for asymptotically almost every $G\sim\mathcal{G}(n,d)$,
\[(1+h_n(\varepsilon))p_{\lambda}\geq \zz{P}_G\{\xi^u_{t_+}\neq \emptyset \}\geq (1-h_n(\varepsilon))p_{\lambda}.\]
\end{Proposition}
\begin{proof}
First of all $\zz{P}_G\{\xi^u_{t_+}\neq \emptyset \}\leq(1+o(1))p_{\lambda}$. This is because
\begin{equation}\label{density.treebound}
\zz{P}_G\{\xi^u_{t_+}\neq \emptyset\}\leq \mathbb{P}\{\tilde{\xi}_{t_+}\neq \emptyset\}=(1+o(1))p_{\lambda},
\end{equation}
where $(\tilde{\xi}_{t})_{t\geq 0}$ is a contact process
on $\zz{T}_{d}$ with the root as the initial configuration.

It remains to show the lower bound. Denote
\[S_u =\sum_
{v\in V(G), v\neq u}\mathbf{1}_{\{v\in \xi^u_{t_+}\}}.\]
In the remains of this section, $o(1)$ terms only depend on $n$
but not $u$ or $G$. Since $u$ is a good vertex,
\[\mathbb{E}_G S_u\geq(1-o(1))p_{\lambda}^2 n.\]
On the other hand,
\[S_u^2=\sum_{v\in [n]\backslash \{u\}}\sum_{w\in [n]\backslash\{u\}}
\mathbf{1}_{\{v\in \xi^u_{t_+},w\in \xi^u_{t_+}\}},\]
and
\[\mathbb{E}_G S_u^2=\sum_{v\in [n]\backslash \{u\}}\sum_{w\in [n]\backslash\{u\}}\zz{P}_G{\{v\in \xi^u_{t_+},w\in\xi^u_{t_+}\}}.\]

Let $t_+=t_{+,1}+t_{+,2}$, where $t_{n,1}= t_{n,2}=t_+/2$,
Using the graphical representation, the event $\{v\in \xi^u_
{t_+}, u\in \xi^u_{t_+}\}$ happens if and
only if there exits open paths starting from $u$ and reach
both $v$ and $w$ in time $t_+$. That requires both of the
following events to happen:
\begin{compactenum}
\item[(1)] $u$ infects some (random) subset $Z\subset [n]$
at time $t_{+,1}$;
\item[(2)] $Z$ infect both $v$ and $w$ in time interval $[t_{+,1},t_+]$.
\end{compactenum}
Let $\xi^u_t$, $\xi^v_t$, $\xi^w_t$ be 3 mutually independent
contact processes, with initial configurations being $\{u\}$,
$\{v\}$ and $\{w\}$ respectively.
By duality we can reverse the time axis of the second event,
and observing the above two events is no easier than observing
the following two events:
\begin{compactenum}
\item[(3)] $\xi^u_t$ survives to time $t_{+,1}$;
\item[(4)] $\xi^v_t$ and $\xi^w_t$ both survive to time $t_{+,2}$.
\end{compactenum}
It is easy to see that the $\zz{P}_G$-probability of observing (3) and (4),
is no bigger than $(1+o(1))p_\lambda^3$, because of (\ref{density.treebound}).

Therefore, for each pair $(v,w)$, we have
\[\zz{P}_G{\{v\in \xi^u_{t_+},w\in\xi^u_{t_+}\}}\leq (1+o(1))p^3_{\lambda},\]
which implies
\[\mathbb{E}_G S^2_u\leq (1+o(1))p_{\lambda}^3n^2.\]
Now since $u$ is a good vertex,
\[\mathbb{E}_G\{S_u\vert S_u>0\}=\frac{\mathbb{E}_G S_u}{\zz{P}_G\{S_u>0\}}\geq \frac{(1-o(1))p_{\lambda}^2n}{\mathbb{P}_G\{S_u>0\}},\]
while
\[\mathbb{E}_G\{S_u^2\vert S_u>0\}=\frac{\mathbb{E}_GS_u^2}{\mathbb{P}_G\{S_u>0\}}\leq \frac{(1+o(1))p_{
\text{\rm survival}}^3n^2}{\mathbb{P}_G\{S_u>0\}}.\]
However by Jensen's inequality, $\mathbb{E}_G\{S_u\vert S_u>0\}^2\leq \mathbb{E}_G\{S_u^2\vert S_u>0\}$, and we must have
\[\mathbb{P}_G\{S_u>0\}\geq (1-o(1))p_{\lambda}.\]
\end{proof}

From the proof of the above proposition we can also obtain
a good estimate of the size of $\xi^u_{t_+}$.
\begin{Proposition}\label{2.typical survival size}
Suppose $u\in [n]$ is a good vertex, and $\xi^u_t$ is a
contact process with initial state $\{u\}$. Fix $\varepsilon>0$.
Then there exist constants $k_n(\varepsilon)\rightarrow 0$ as
$n\rightarrow \infty$ such that  for asymptotically almost
every $G$,
\[ \zz{P}_G\{(1-\delta)np_{\lambda}\leq S_u\leq
(1+\delta)np_{\lambda}|S_u>0 \}\geq 1-k_n(\varepsilon).\]
\end{Proposition}
\begin{proof}
From Proposition \ref{2.typical survival},
$\mathbb{P}_G \{S_u>0\}\geq (1-o(1))p_{\lambda}$.
On the other hand $\mathbb{P}_G \{S_u>0\}\leq (1+o(1))p_{\lambda}$.
These two bounds combined with the calculation in the proof of
Proposition \ref{2.typical survival}, we obtain
\[\mathbb{E}_G\{S_u^2\vert S_u>0\}\leq(1+o(1))\mathbb{E}_G\{S_u\vert S_u>0\}^2.\]
Then we apply Chebyshev's inequality.
\end{proof}

Now we are ready to prove Theorem \ref{THM2}.
\begin{proof}[Proof of Theorem \ref{THM2}]
Let $S=\sum_{v\in [n]}\mathbf{1}_{\{v\in \xi^G_{t_+}\}}$. By duality of
the contact process, Proposition \ref{2.typical}
and Proposition \ref{2.typical survival},
\[\mathbb{E}_GS\geq (1-o(1))np_{\lambda}.\]
On the other hand, using a similar argument as
in the proof of  Proposition \ref{2.typical survival}, we have
\[\mathbb{E}_GS^2\leq (1+o(1))n^2p_{\lambda}^2.\]
Therefore we obtain
\[\mathbb{E}_GS^2\leq(1+o(1))(\mathbb{E}_G S)^2.\]
Then we apply Chebyshev's inequality.
\end{proof}
\section*{Acknowledgements}
The authors would like to thank Si Tang for her help in drawing the figures.

\bibliographystyle{plain}
\bibliography{mainbib}

\end{document}